\documentclass[12pt]{amsart}
\usepackage{amssymb,latexsym,eufrak,amsmath,amscd,
graphicx}

\theoremstyle{plain}
\newtheorem{theorem}{Theorem}[section]
\newtheorem{proposition}[theorem]{Proposition}
\newtheorem{corollary}[theorem]{Corollary}
\newtheorem{lemma}[theorem]{Lemma}

\newtheorem{theoremalpha}{Theorem}

\theoremstyle{definition}
\newtheorem{definition}[theorem]{Definition}
\newtheorem{remark}[theorem]{Remark}
\newtheorem{example}[theorem]{Example}

\newcommand{\ra}{\rightarrow}
\newcommand{\lra}{\longrightarrow}
\newcommand{\noi}{\noindent}
\newcommand{\PP}{\mathbf{P}}

\newcommand{\RR}{\mathbf{R}}
\newcommand{\NN}{\mathbf{N}}

\newcommand{\ZZ}{\mathbf{Z}}
\newcommand{\CC}{\mathbf{C}}
\newcommand{\QQ}{\mathbf{Q}}

\newcommand{\OO}{\mathcal  {O}}

\newcommand{\ci}{\mathcal  {I}}

\newcommand{\fra}{\frak{a}}
\newcommand{\tn}[1]{\textnormal{{#1}}}
\newcommand{\frb}{\frak{b}}

\newcommand{\frmm}{\frak{m}}

\newcommand{\linser}[1]{\vert \mspace{1.75mu} {#1}
\mspace{1.75mu} \vert}
\newcommand{\alinser}[1]{\Vert\mspace{1.5mu} {#1}
\mspace{1.5mu} \Vert}

\newcommand{\HH}[3]{H^{{#1}} \big( {#2} , {#3}
\big) }
\newcommand{\hh}[3]{h^{{#1}} \big( {#2} , {#3}
\big) }

\newcommand{\MI}[1]{\mathcal  {J} ( {#1} ) }

\newcommand{\AMI}[1]{\mathcal{J}\big( \mspace{1.75mu}
                            \alinser{#1} \mspace{1.75mu} \big)}
\newcommand{\AMMI}[2]{\MI{{#1}, \alinser{{#2}}}}

\newcommand{\pr}{\prime}

\newcommand{\Div}{\text{Div}}

\DeclareMathOperator{\vol}{vol}

\DeclareMathOperator{\ord}{ord}

\DeclareMathOperator{\BB}{{\bf B}}

\DeclareMathOperator{\B+}{\BB_{+}}

\newcommand{\bs}[1]{\frb(  {#1} )}

\begin{document}

\title{Restricted volumes and  base loci of linear
series}

\author[L.Ein]{Lawrence Ein}
\address{Department of Mathematics \\ University of
Illinois at Chicago, \hfil\break\indent  851 South
Morgan Street (M/C 249)\\ Chicago, IL 60607-7045, USA}
\email{ein@math.uic.edu}

\author[R. Lazarsfeld]{Robert Lazarsfeld}
\address{Department of Mathematics
\\ University of Michigan \\ Ann Arbor, MI 48109, USA}
\email{rlaz@umich.edu}

\author[M. Musta\c{t}\u{a}]{Mircea Musta\c{t}\u{a}}
\address{Department of Mathematics
\\ University of Michigan \\ Ann Arbor, MI 48109, USA}
\email{mustata@umich.edu}

\author[M. Nakamaye]{Michael Nakamaye}
\address{Department of Mathematics and Statistics\\
University of New
  Mexico,
\hfil\break\indent Albuquerque, New Mexico 87131, USA}
\email{nakamaye@math.unm.edu}

\author[M. Popa]{Mihnea Popa}
\address{Department of Mathematics \\
  University of Chicago \\ 5734 S. University Ave.,
\hfil\break\indent Chicago, IL 60637,
  USA}
\email{mpopa@math.uchicago.edu}

\begin{abstract} We introduce and study the restricted
volume of a divisor along a subvariety. Our main result
is a description of the irreducible components of the
augmented base locus by the vanishing of the restricted
volume.
\end{abstract}

\thanks{2000\,\emph{Mathematics Subject Classification}.
 Primary 14C20; Secondary 14E15, 14E05.
\newline The research of Ein, Lazarsfeld,
Musta\c{t}\u{a} and Popa
  was partially supported by the NSF under grants DMS
0200278, DMS 0652845, DMS 0500127, and DMS 0200150}
\keywords{Volume, Seshadri constant, base locus, big
divisor}

\maketitle

\section*{Introduction}

Let $X$ be a smooth complex projective variety of dimension $n$.
While it is classical that ample line bundles on $X$ display
beautiful geometric, cohomological and numerical properties, it was
long believed that one couldn't hope to say much in general  about
the behavior of arbitrary effective divisors. However, it has
recently become clear (\cite{nakamaye1}, \cite{nakayama},
\cite{positivity}, \cite{elmnp1}),  that many aspects of the
classical picture do in fact extend to arbitrary effective (or
``big" divisors) provided that one works asymptotically.  For
example, consider the \textit{volume} of a divisor $D$:
\[ \vol_X(D) \ =_\text{def} \ \limsup_{m\to\infty}  \ \frac{
\hh{0}{X}{\OO_X(mD)}}{m^n/n!}. \]
 When $A$ is ample, it follows from the
asymptotic Riemann-Roch formula
 that the volume is just the top
self-intersection number of $A$:
\[  \vol_X(A) \ = \ \big( A^n \big). \]
In general, one can view $\vol_X(D)$ as the natural generalization
to arbitrary divisors of this self-intersection number. (If
$D$ is not ample, then the actual intersection number $( D^n)$
typically doesn't carry immediately useful geometric information.
For example, already on surfaces it can happen that $D$ moves in a
large linear series while $(D^2) \ll 0.$) It turns out that many of
the classical properties of the self-intersection number for ample
divisors extend in a natural way to the volume. For instance, it was
established by the second author in \cite{positivity} that
$\vol_X(D)$ depends only on the numerical equivalence class of $D$,
and that it determines a continuous function
\[ \vol_X : N^1(X)_{\RR} \lra \RR \]
on the finite dimensional vector space of numerical
equivalence classes of $\RR$-divisors.

Now consider an irreducible subvariety $V \subseteq X$
of dimension $d$. In the classical setting, when $A$ is
ample, the intersection numbers $(A^d \cdot V)$ play an
important role in many geometric questions. The goal of
the present paper is to study the asymptotic analogue
of this degree for an arbitrary divisor $D$.
Specifically, the \textit{restricted volume} of $D$
along
$V$ is defined to be
$$\vol_{X\vert V}(D)\ = \
\limsup_{m\to\infty}\frac{\dim_{\CC}\, {\rm
Im}\Big(\HH{0}{X}{\OO_X(mD)} \lra
\HH{0}{V}{\OO_V(mD)}\,  \Big)} {m^d/d!}.$$ Thus
$\vol_{X|V}(D)$ measures asymptotically the number of
sections of the restriction
$\OO_V(mD)$ that can be lifted to $X$.  For example, if
$A$ is ample then the restriction maps are eventually
surjective, and hence
\[  \vol_{X|V}(A) \ = \ \vol_V(\,   A_{|V}\, ) \ = \
\big( A^d \cdot V \big). \] In general however it can happen that
$\vol_{X|V}(D) < \vol_V(D_{|V})$. The definition extends in the
evident way to $\QQ$-divisors. Restricted volumes seem to have first
appeared in passing in Tsuji's preprint \cite{tsuji1}, and they play
an important role in the papers \cite{hm}, \cite{takayama2} of
Hacon--McKernan and Takayama elaborating Tsuji's work.

In order to state our results, we need to be able to
discuss  how
$V$ sits with respect to the base-loci of $D$.
Recall to this end  that the
\textit{stable base-locus} $\BB(D)$ of an integral or
$\QQ$-divisor
$D$ is by definition the common base-locus of the
linear series $\linser{mD}$ for all sufficiently large
and divisible $m$. Unfortunately, these loci behave
rather unpredictably: for example, they don't depend
in general only on the numerical equivalence class of
$D$. The fourth author observed in \cite{nakamaye1} that one
obtains a much cleaner picture if one perturbs $D$
slightly by
 subtracting off a small ample class.  Specifically,
the
\textit{augmented base-locus} of $D$ is defined to be
\[
\B+(D) \ =_{\text{def}} \ \BB(D-A) \]
for a small ample $\QQ$-divisor $A$, this being
independent of $A$ as long as its class in
$N^1(X)_{\RR}$ is sufficiently small. Thus
 $\B+(D)
\supseteq
\BB(D)$. These augmented base-loci were studied
systematically in \cite{elmnp}, where in particular it was
established that $\B+(D)$ depends only on the numerical
equivalence class of $D$. Since the
definition involves a perturbation,
$\B+(\xi)$ is consequently also defined for any  class
$\xi
\in N^1(X)_{\RR}$.

Our first result  involves the formal behavior of the
restricted volume
$\vol_{X|V}(D)$.
\begin{theoremalpha}\label{thmA} Let $V \subseteq X$ be
an irreducible subvariety of dimension $d > 0$ and let
$D$ be a $\QQ$-divisor such that $V \not \subseteq
\B+(D)$.  Then \[ \vol_{X|V}(D) > 0, \] and
$\vol_{X|V}(D)$ depends only on the numerical
equivalence class of $D$.  Furthermore, $
\vol_{X|V}(D)$ varies continuously as a function of the
numerical equivalence class of $D$, and it extends
uniquely to a continuous function
\[  \vol_{X|V} : \textnormal{Big}^V(X)^+_\RR \lra \RR,
\] where $\textnormal{Big}^V(X)^+_\RR$ denotes the set
of all real divisor classes $\xi$ such that $V \not \subseteq
\B+(\xi)$. This function is homogeneous of degree $d$, and it
satisfies the log-concavity property
\[
\vol_{X|V}(\xi_1 + \xi_2)^{1/d} \ \ge \
\vol_{X|V}(\xi_1)^{1/d} \, + \, \vol_{X|V}(\xi_2)^{1/d}.
\]
\end{theoremalpha}

We also show that one can compute $\vol_{X|V}(D)$ in
terms of ``moving intersection numbers" of divisors with
$V$:
\begin{theoremalpha}\label{thmB} Assume as above that
$D$ is a
$\QQ$-divisor on $X$, and that $V$ is a subvariety of
dimension
$d > 0$
 such that
$V \not\subseteq \B+(D)$. For every large
and sufficiently divisible integer $m$, choose $d$
general divisors
$E_{m,1}, \ldots, E_{m,d} \in
\linser{mD}$. Then
$$\vol_{X\vert V}(D) \ = \ \lim_{m\to\infty}
\frac{\#\big(\  V\cap E_{m,1}\cap\ldots\cap
E_{m,d} \, \smallsetminus\, \BB(D)\  \big)}{m^d}.$$
\end{theoremalpha}
\noi In other words, $\vol_{X|V}(D)$ computes the rate of growth of
the number of intersection points away from $\BB(D)$  of $d$
divisors in $\linser{mD}$ with $V$. If $D$ is ample, this just
restates the fact that $\vol_{X|V}$ is given by an intersection
number. The theorem  extends one of the basic properties of
$\vol_X(D)$, essentially due to Fujita; as in the case $V = X$, the
crucial point is to show that one can approximate $\vol_{X|V}(D)$
arbitrary closely by intersection numbers with ample divisors on a
modification of $X$  (cf. \cite{positivity} \S11.4.A or
\cite{del}). This result has been proved independently by Demailly
and Takayama \cite{takayama2}. It also leads to an extension of
the theorem of Angehrn and Siu \cite{AS} on effective base-point
freeness of adjoint bundles in terms of  restricted volumes (see
Theorem~\ref{Siu}).

Our main result is  that these restricted volumes
actually govern base-loci. By way of background,
suppose that $P$ is a nef divisor on $X$. The fourth
author proved in \cite{nakamaye1} that the irreducible
components of $\B+(P)$ consist precisely of those
subvarieties
$V$ on which $P$ has degree zero, i.e.
\begin{equation}\label{augmented_nef}
\B+(P) \ = \  \underset{(P^{\dim\,V}\cdot
V)=0}{\bigcup}~V,
\end{equation} where $V$ is required to be positive
dimensional. We prove the analogous result for
arbitrary $\QQ$-divisors
$D$:
\begin{theoremalpha}\label{thmC} If $D$ is a
$\QQ$-divisor on $X$, then $\B+(D)$ is the union of all
positive dimensional subvarieties $V$ such that
$\vol_{X\vert V}(D)=0$.
\end{theoremalpha}
\noi One can extend the statement to $\RR$-divisors by
introducing the set $\textnormal{Big}^V(X)_\RR$
consisting  of all real divisor classes such that $V$ is
not properly contained in any irreducible component of
$\B+(\xi)$. Then
$\vol_{X|V}$ determines  a continuous function
\[  \vol_{X|V} : \textnormal{Big}^V(X)_\RR\lra \RR \]
with the property that $$ \vol_{X\vert V}(\xi) = 0 \iff  V ~{\rm is
  ~an~irreducible~component~of}~\B+(\xi).$$
(Note that just as it can happen for a nef divisor $P$
that $( P^ {\dim V} \cdot V ) = 0$ while $ (
P^{\dim W} \cdot W ) > 0$ for some $W \subseteq V$,
so it can happen that $\vol_{X|V}(D) = 0$ but
$\vol_{X|W}(D) > 0$ for some $W \subseteq V$. This is
why one has to focus here on irreducible components of
base loci. Compare also Example \ref{failure}.)

  The proof of Theorem~\ref{thmC} is based on
ideas introduced by the fourth author in
\cite{nakamaye2}, together with a result
(Theorem~\ref{separation_of_jets}) describing
$\vol_{X\vert V}(D)$ in terms of separation of jets at
general points of $V$. This allows one to lift sections
of line bundles from a subvariety to $X$ as a result of
direct computation rather than vanishing of cohomology.
The very rough idea of the proof is the following:
starting with a lower bound for $\vol_{X\vert V}(D+A)$
for some ample divisor $A$, we deduce that there are
points on
$V$ at which the line bundles
$\OO(m(D+A))$ separate many jets, for
large enough
$m$. This allows us in turn to produce lower bounds for
the dimension of spaces of sections with small
vanishing order at the same points, for line bundles of
the form $\OO(m(D- A'))$, with $A'$ a new ample
divisor. The conclusion is that the asymptotic
vanishing order of
$D-A'$ along $V$ (${\rm
 ord}_V(\Vert
D-A'\Vert)$ in the notation of \cite{elmnp}) can be
made very small. However, as we make
$A \rightarrow 0$, we prove that there exists a uniform
constant $\beta > 0$ such that if
$V\subseteq \B+(D)$, then ${\rm ord}_V(\Vert D-A'\Vert)
>
\beta\cdot\Vert A'\Vert$, which produces a contradiction. The actual
proof is quite technical and occupies most of \S 5.

In the last section we make the connection between our
results and those of \cite{nakamaye2}, describing the
augmented base locus in terms of another asymptotic
invariant, the \emph{moving Seshadri constant}. This
invariant was introduced in \cite{nakamaye2} as a
generalization to arbitrary big divisors of the usual
notion of Seshadri constant for big and nef divisors
(cf. \cite{positivity} \S 5.1) . We describe the
relationship between moving Seshadri constants and
restricted volumes, and as a consequence of
the results  above we obtain a slight
strengthening of the main result in \cite{nakamaye2}:
\emph{the moving Seshadri constant varies as a
continuous function on $N^1(X)_{\RR}$, and given an
arbitrary $\RR$-divisor $D$,
$\B+(D)$ is the set of points at which the moving
Seshadri constant of $D$ is zero}.

The results in this paper are part of a more general program of
using asymptotic invariants of divisors in order to get information
about the geometry of linear series, base loci, and cones of
divisors on a projective variety. Invariants of a different flavor
were used in \cite{elmnp} in order to describe a lower approximation
of the stable base locus of a divisor, called the \emph{restricted
base locus} (or \emph{non-nef locus}, cf. \cite{boucksom1},
\cite{bdpp}, see also \cite{debarre}). The reader can also find
there a thorough discussion of the connections between these various
asymptotic base-locus-type constructions. Finally, we refer to
\cite{elmnp1} for an overview of the basic ideas revolving around
asymptotic invariants of line bundles.

\section{The augmented base locus}

We start by fixing some notation. Let $X$ be a smooth complex
projective variety of dimension $n$. An integral divisor $D$ on $X$
is an element of the  group $\Div(X)$ of Cartier divisors. The
corresponding linear series is denoted by
 $|D|$ and its base locus by ${\rm Bs}(|D|)$. As usual
we can speak about $\QQ$- or $\RR$-divisors. A
$\QQ$- or
$\RR$-divisor
$D$ is
\textit{effective} if it is a non-negative linear
combination of effective integral divisors with $\QQ$-
or $\RR$-coefficients.
  If $D$ is effective, we denote by
${\rm Supp}(D)$ the union of the irreducible components
which appear in the associated
 Weil divisor. We often use the same notation for an
integral divisor and for the corresponding line bundle.
Numerical equivalence between $\QQ$- or
$\RR$-divisors will be denoted by $\equiv$. We denote
by $N^1(X)_{\QQ}$ and
$N^1(X)_{\RR}$ the finite dimensional $\QQ$- and
$\RR$-vector spaces of numerical equivalence classes.
One has $N^1(X)_{\RR} = N^1(X)_{\QQ}
\otimes_{\QQ} \RR$, and we fix compatible norms
$\parallel\cdot\parallel$ on these two spaces. Given a
divisor $D$, we write $\Vert D \Vert$ for the norm of
the class of $D$.  A $\QQ$-divisor is \emph{big} if for
$m$ divisible enough, the linear series $|mD|$ defines
a birational map onto its image. One can show that
$D$ is big if and only if $D\equiv A+E$, where $A$ is
ample and $E$ is effective. This can be taken as
definition in the case of an $\RR$-divisor (see
\cite[Section 2.2]{positivity}  for the basic
properties of big divisors). The
\emph{big cone} is the open convex cone in
$N^1(X)_{\RR}$ consisting of big $\RR$-divisor classes.

We recall from \cite{elmnp} the definition of the
augmented base locus of a divisor. Suppose first that
$D$ is a $\QQ$-divisor on $X$. The \emph{stable base
locus} of $D$ is
$$\BB(D):=\bigcap_{m}{\rm Bs}(|mD|)_{\rm red},$$ where
the intersection is over all $m$ such that $mD$ is an
integral divisor. It is easy to see that if $p$ is
divisible enough, then
$\BB(D)={\rm Bs}(|pD|)_{\rm red}$.

The \emph{augmented base locus} of an $\RR$-divisor $D$
is defined to be
$$\B+(D):=\bigcap_A\BB
(D-A),$$ where the intersection
is over all ample divisors $A$ such that $D-A$ is a
$\QQ$-divisor. Equivalently, we have
$$\B+(D)=\bigcap_{D=A+E}{\rm Supp}(E),$$ where we take
the intersection over all decompositions $D=A+E$, with
$A$ ample and $E$ effective. It follows from definition
that $\B+(D)$ is a closed subset of $X$, and
$\B+(D)\neq X$ if and only if $D$ is big. Moreover,
there is $\eta>0$ such that for every ample divisor $A$
with $\parallel A\parallel<\eta$ and such that $D-A$ is
a $\QQ$-divisor, we have $\B+(D)=\BB(D-A)$. It is easy
to see that if $D\equiv E$ then $\B+(D)=\B+(E)$. For a
detailed study of augmented base loci, see \cite{elmnp}
\S1. In addition, we will need the following property.

\begin{proposition}\label{no_points} If $D$ is a
$\QQ$-divisor, then $\BB(D)$ has no isolated points. In particular,
for every $\RR$-divisor $D$, the augmented base locus $\B+(D)$ has
no isolated points.
\end{proposition}

\begin{proof}  Suppose that $x$ is an isolated
point in $\BB(D)$, and let $m$ be large and divisible enough such
that $mD$ is integral and $\BB(D)={\rm Bs}(|mD|)_{\rm red}$. Let
$\fra\subseteq {\mathcal O}_X$ be the ideal defining the scheme
${\rm Bs}(|mD|)\smallsetminus \{x\}$ and let $f\colon X'\to X$ be
the normalized blow-up along $\fra$. We can write $f^*(mD)=M+F$,
where $f^{-1}(\fra)={\mathcal O}(-F)$ and the base locus of $|M|$ is
concentrated at the point $f^{-1}(x)$. A result of Zariski (see
\cite{zariski}, and also \cite{ein}) implies that there is $p$ such
that $|pM|$ is base-point free. Therefore $x$ is not in the base
locus of $|pmD|$, a contradiction.

If $D$ is an $\RR$-divisor, then $\B+(D)=\BB(D-A)$ for some ample
divisor $A$ such that $D-A$ is a $\QQ$-divisor, so the last
assertion follows.
\end{proof}

\begin{remark}
The assertion on $\B+(D)$ in the above proposition can also be
proved using the elementary theory of multiplier ideals,
 avoiding the appeal to Zariski's theorem.
\end{remark}

\section{Restricted volumes and asymptotic intersection
numbers}

\noindent {\bf Restricted volumes.} Recall that $X$ is a smooth
projective variety. For any line bundle $L$ on $X$ and any
subvariety $V\subseteq X$, we set
$$\HH{0}{X|V}{L} \ := \ {\rm
Im}\Big(\HH{0}{X}{L} \longrightarrow \HH{0}{V}{
L|_{V}}\Big),$$ while
$\hh{0}{X|V}{L}$ is the dimension of
$\HH{0}{X|V}{L}$. Of course we use the analogous
notation for divisors.

\begin{definition}[Restricted volume] If $L$ is a line
bundle on $X$, and if $V\subseteq X$ is a subvariety of
dimension $d\geq 1$, then the \emph{restricted volume}
of $L$ along $V$ is
\begin{equation}\label{definition_volume}
\vol_{X\vert V}(L)\:=\ \limsup_{m\to\infty}\frac{
\hh{0}{X|V}{mL}}{m^d/d!}.
\end{equation}
\end{definition}
\noi Again, the same definition applies to divisors.

Note that if $V=X$, then the restricted volume of $L$
along $V$ is the usual volume of $L$, denoted by
$\vol_X(L)$, or simply by $\vol(L)$. We refer to
\cite{positivity} \S 2.2.C
 for a study of the volume function. Our main
goal in this section is to extend these results to the
case of an arbitrary subvariety $V\subseteq X$. As we
will see, everything goes over provided  we
assume that $V\not\subseteq\B+(L)$.  To begin with, the
following lemma implies immediately that
$$\vol_{X\vert V}(qL)=q^d\vol_{X\vert V}(L),$$ so we
can also define in the obvious way $\vol_{X\vert V}(D)$
when $D$ is a $\QQ$-divisor.

\begin{lemma}\label{scaling} Let $D$ be any divisor
on $X$ and $q\in \NN$ a fixed positive integer. Then
$$\limsup_{m\to\infty}\frac{\hh{0}{X|V}{mD}}{m^d/d!} =
\limsup_{m\to\infty}\frac{\hh{0}{X|V}{qmD}}{(qm)^d/d!}.$$
\end{lemma}
\begin{proof} The proof is identical to that of the
corresponding statement for the usual volume function
of a line bundle given in \cite{positivity} Lemma
2.2.38.
\end{proof}

\begin{example}[Ample and nef divisors]\label{ample} If
$D$ is an ample $\QQ$-divisor, then Serre vanishing
implies
$\vol_{X\vert V}(D)=\vol_V(D\vert_V)=(D^d\cdot V)$. We
will see later that the same thing is true if $D$ is
only nef under the hypothesis
$V\not\subseteq \B+(D)$ (cf. Corollary
\ref{volume_for_nef} and Example~\ref{big_and_nef}).
\end{example}

\begin{lemma}\label{pull_back} Let $f :
X'\longrightarrow X$ be a proper, birational morphism of
smooth varieties and let $D$ be a $\QQ$-divisor on $X$.
If $V'\subseteq X'$ and $V=f(V')$ have the same
dimension, then \[ \vol_{X'\vert
V'}(f^*(D))\ = \ \vol_{X\vert V}(D). \]
\end{lemma}

\begin{proof} It is enough to note that for every $m$
such that $mD$ is an integral divisor, we have the
commutative diagram
\[
\begin{CD} H^0(X,mD)@>u>>H^0(X',mf^*(D))\\ @VVV  @VVV\\
H^0(V,mD\vert_V)@>v>>H^0(V',mf^*(D)\vert_{V'}),
\end{CD}
\] where $u$ is an isomorphism and $v$ is a
monomorphism.
\end{proof}

\begin{remark} Note by contrast that
\[ \vol_{V^\pr}( f^*(D)|_{V^\pr}) \ = \ \deg(V^\pr
\lra V) \cdot \vol_V(D|_V). \]
\end{remark}

\noindent {\bf Asymptotic intersection numbers.} Let
$D$ be a $\QQ$-divisor. Note that if $V\subseteq
\BB(D)$, then clearly $\vol_{X\vert V}(D) = 0$. Assume
now that $V\not\subseteq \BB(D)$. We can then define
another invariant, an asymptotic intersection number of
$D$ and $V$, in the following way. Fix a natural number
$m > 0$ which is sufficiently divisible so that
$\BB(D)={\rm Bs}(|mD|)_{\rm red}$, and let
\[ \pi_m: X_m \ra X\] be a resolution of the base ideal
$\frb_m = \bs{\linser{mD}}$. Thus we have a
decomposition
\[ \pi_m^\ast \big(\linser{mD} \big) \ = \ \linser{M_m}
+ E_m,
\] where
${M_m}$ (the \textit{moving part} of $\linser{mD}$) is
free, and
$E_m$ is the \textit{fixed part}. We can  --- and
without further mention, always will ---  choose all
such resolutions with the property that they are
isomorphisms over the generic point of
$V$. We then denote by $\widetilde{V}_m$ the proper
transform of $V$, which by hypothesis is not contained
in ${\rm Supp}(E_m)$.
\begin{definition}[Asymptotic intersection number]
\label{Asymptotic intersection number} With the
notation just introduced, the
\emph{asymptotic intersection number} of
$D$ and $V$ is defined to be
$$\alinser{ D^d \cdot V} \  := \
\limsup_{m\to\infty}\frac{(M_m^d\cdot
\widetilde{V}_m)}{m^d}.$$
\end{definition}
\noi Naturally enough, we make the analogous definition
for line bundles.

\begin{remark} The intersection numbers with   $M_m$
have the following interpretation: if
$D_1,\ldots,D_d$ are general divisors in $|mD|$, then
$(M_m^d\cdot \widetilde{V}_m)$ is equal to the number
of points in $D_1\cap\ldots\cap D_d\cap V$ that do not
lie in ${\rm Bs}(|mD|)$. In particular, this number
does not depend on the resolution we are choosing.
\end{remark}

\begin{remark} Another sort of asymptotic intersection
number, involving the moving intersection points of $n$ different
big line bundles on $X$, is introduced and studied in \cite{bdpp}.
Under suitable conditions on the position of $V$ with respect to the
relevant base-loci, one could combine the two lines of thought to
define an asymptotic intersection number of the type
\[  \alinser{L_1 \cdot   \ldots \cdot L_d \cdot V}. \]
However we do not pursue this here.
\end{remark}

\begin{remark}
\label{Mov.Int.No.Is.Limit} We note that
$\alinser{ D^d
\cdot V}$ computes in fact the \emph{limit}
$\lim_{m\to\infty}\frac{(M_m^d\cdot
\widetilde{V}_m)}{m^d}=\sup_m\frac{(M_m^d\cdot\widetilde{V}_m)}{m^d}$,
where the limit and the supremum are over all
$m$ such that $\BB(D)={\rm Bs}(|mD|)_{\rm red}$.
Indeed, given such $p$ and $q$, we may take
$\pi : X'\longrightarrow X$ that satisfies our
requirements for
$|pD|$, $|qD|$ and $|(p+q)D|$. If $\widetilde{V}$ is
the proper transform of $V$ on $X'$, we deduce that
$M_{p+q}-(M_p+M_q)$ is effective and does not contain
$\widetilde{V}$ in its support, so
$$(M_{p+q}^d\cdot\widetilde{V})^{1/d}\ \geq \
((M_p+M_q)^d\cdot\widetilde{V})
^{1/d}\
\geq
\
(M_p^d\cdot\widetilde{V})^{1/d}+(M_q^d\cdot\widetilde{V})^{1/d},$$
where for the last inequality we refer to
\cite{positivity}, Corollary 1.6.3. It is standard to
deduce our claim from this inequality.
\end{remark}

\begin{remark}It follows from the previous remark that
$\alinser{ (mD)^d\cdot V} = m^d\parallel
D^d\cdot V\parallel$ for every $m$.
\end{remark}

The next result gives another interpretation
of these
intersection numbers.
\begin{proposition}\label{general_resolution} If $D$ is
a $\QQ$-divisor, and if $V$ is not contained in
$\B+(D)$, then
\begin{equation}\label{formula_for_volume}
\alinser{ D^d \cdot V} \
= \\sup_{\pi^*D=A+E}(A^d\cdot{\widetilde{V}}),
\end{equation} where the supremum is taken over all
projective birational morphisms
$\pi : X'\longrightarrow X$ with $X'$ smooth, that give
an isomorphism at the general point of $V$, and over
all expressions $\pi^*D=A+E$, where $A$ and $E$ are
$\QQ$-divisors, with $A$ ample, $E$ effective and
$\widetilde{V}\not\subseteq {\rm Supp}(E)$. {\rm (}Here
$\widetilde{V}$ denotes the proper transform of
$V$.{\rm )}
\end{proposition}

\begin{proof} Consider first any morphism $\pi$ as in
the statement of the proposition and let $m$ be
divisible enough. The number of points outside ${\rm
Bs}(|mD|)$ that lie on the intersection of $d$ general
members of $|mD|$ with $V$ is the same as the number of
points outside ${\rm Bs}(|\pi^*(mD)|)=\pi^{-1}({\rm
Bs}(|mD|))$ that lie on the intersection of
$d$ general members of $|m\pi^*(D)|$ with
$\widetilde{V}$. Moreover, since $\widetilde{V}$ is not
contained in ${\rm Supp}(E)$, this number is at least
the number of intersection points of $\widetilde{V}$
with
$d$ general members in $|mA|$, which is
$m^d(A^d\cdot\widetilde{V})$. Dividing by $m^d$ and
letting $m$ go to infinity gives the inequality
``$\geq$" in the statement.  On the other hand, by
definition we have
$$\alinser{D^d \cdot V}   \ := \
\limsup_{m\to\infty}\frac{\big( M_m^d\cdot
\widetilde{V}_m\big)}{m^d}.$$ It is easy to see that
since
$V\not\subseteq\B+(D)$, we have
$\widetilde{V}_m\not\subseteq\B+(M_m)$. Therefore we
can write $M_m=A+E$, with $A$ ample, $E$ effective, and
$\widetilde{V}_m\not\subseteq{\rm Supp}(E)$. If
$p\in\NN^*$, then we have
$M_m=(1/p)E+A_p$, where
$A_p=\frac{1}{p}A+\frac{p-1}{p}M_m$ is ample since
$M_m$ is nef. The opposite inequality in the statement
follows from
$\lim_{p\to\infty}(A_p^d\cdot\widetilde{V}_m)=(M_m^d\cdot\widetilde{V}_m)$.
\end{proof}

As the right hand side of (\ref{formula_for_volume})
depends only on the numerical class of $D$, we deduce
the following:
\begin{corollary}\label{invariance_asymp_int} If
$D_1\equiv D_2$ are $\QQ$-divisors, and if
$V\not\subseteq\B+(D_1)$, then
$$\alinser{D_1^d\cdot V} \ = \ \alinser{ D_2^d\cdot
V}.$$
\end{corollary}

\noindent {\bf A generalized Fujita Approximation
Theorem.} The next result shows that if $V$ is not
contained in $\B+(D)$, then the two invariants we have
defined for $D$ along $V$ are the same. In the case
$V=X$, this is Fujita's Approximation Theorem (see
\cite{del}). In addition, we give a formula for the
restricted volume in terms of asymptotic multiplier
ideals, connecting our approach to ideas for defining
invariants due to Tsuji \cite{tsuji2}
We mention that the relationship between asymptotic
intersection numbers
  and the growth of sections vanishing along restricted
multiplier
  ideals appears also in the recent work of
Takayama \cite{takayama2}. Note that these statements
are interesting only for big divisors
$D$, since otherwise $\B+(D)  = X$. We will assume
familiarity with the basic theory of multiplier ideals
developed in Part III of \cite{positivity}.

If $D$ is an integral divisor, we denote by
$\AMI{mD}$ the asymptotic multiplier
ideal of $mD$.  For
simplicity, we use the following notation: if $\ci$ is
an ideal sheaf on a variety $X$, and if $V\subseteq X$
is a subvariety, then  $\ci\vert_V$ denotes the ideal
$\ci\cdot\OO_V$.

\begin{theorem}\label{generalized_Fujita} Let  $D$ be a
$\QQ$-divisor on $X$, and let $V$ be a $d$-dimensional subvariety of
$X$ \tn{(}$d \ge 1\tn{)}$ such that $V\not\subseteq\B+ (D)$.  Then
$$\vol_{X\vert V}(D)\ = \  \alinser{D^d\cdot V} \ = \
\limsup_{m\to\infty}\frac{\hh{0}{V}{ \OO(mD)\otimes
\AMI{mD}\vert_V}}{m^d/d!},$$ where in the last term we take the
limit over $m$ sufficiently divisible so that $mD$ is integral.
\end{theorem}

The proof of the above theorem will be given in the next section. We
record now several consequences and examples.
Theorem~\ref{generalized_Fujita} together with
Corollary~\ref{invariance_asymp_int} imply the following:
\begin{corollary}\label{numerical_invariance} If $D$ is
a $\QQ$-divisor and if
$V$ is a subvariety such that $V\not\subseteq\B+(D)$,
then the restricted volume $\vol_{X\vert V}(D)$ depends
only on the numerical class of $D$.
\end{corollary}

\begin{corollary}\label{lim} If $D$ is a $\QQ$-divisor
on
$X$ and $V\subseteq X$ is a $d$-dimensional subvariety
such that
$V\not\subseteq \B+ (D)$, then
$$\limsup_{m\to\infty}\frac{\hh{0}{X|V}{mD}}{m^d/d!} \ =
\  \lim_{m\to\infty}\frac{\hh{0}{X|V}{mD}}{m^d/d!}.$$ In other
words, the restricted volume is actually the limit
$$\vol_{X\vert V}(D) =
\lim_{m\to\infty}\frac{\hh{0}{X|V}{mD}}{m^d/d!}.$$
\end{corollary}
\begin{proof} This fact is an immediate consequence of
Theorem~\ref{generalized_Fujita} together with
Proposition~\ref{general_resolution}, the proof being
identical to that of the corresponding statement for
the usual volume (\cite{positivity}, Example~11.4.7).
\end{proof}

\begin{corollary}\label{concavity} Let $D_1$ and $D_2$
be two $\QQ$-divisors and $V\subseteq X$ a subvariety
of dimension $d\geq 1$ such that
$V\not\subseteq\B+(D_1)\cup\B+(D_2)$. Then
$$\vol_{X\vert V}(D_1+D_2)^{1/d}\ \geq\ \vol_{X\vert V}
(D_1)^{1/d}\, + \, vol_{X\vert V}(D_2)^{1/d}.$$
\end{corollary}

\begin{proof} The assertion follows from
Theorem~\ref{generalized_Fujita} and
Proposition~\ref{general_resolution}, using the
corresponding concavity property for the volumes of
ample line bundles (see \cite{positivity}, Corollary
1.6.3).
\end{proof}

\begin{corollary}\label{volume_for_nef} Suppose that
$D$ is a nef $\QQ$-divisor and that $V\subseteq X$ is a
subvariety of dimension $d\geq 1$. If $V$ is not
contained in $\B+(D)$, then $\vol_{X\vert
V}(D)=(D^d\cdot V)$.
\end{corollary}

\begin{proof} We may assume that $D$ is an integral
divisor. Since $D$ is nef and big, we have
$$\AMMI{X}{mD} \ = \ \OO_X~ {\rm for~ all}~
m$$ (\cite{positivity}, Proposition 11.2.18). It
follows from Theorem~\ref{generalized_Fujita} that
$$\vol_{X\vert V}(D)\ =\ \vol_V(D\vert_V)\ = \ (D^d\cdot
V),$$ where we use the corresponding result for the
usual volume function (\cite{positivity},
Corollary~1.4.41).
\end{proof}

\begin{remark}\label{remark_nef} If $D$ is a nef
divisor, then $\B+(D)={\rm Null}(D)$, where
${\rm Null}(D)$ is the union of the subvariaties $V$ of
$X$ such that
$D\vert_V$ is not big. This is the main result of
\cite{nakamaye1}, which we will reprove in
Corollary~\ref{cor_nakamaye} below, allowing also
$\RR$-coefficients.
\end{remark}

\begin{example}\label{Zariski} Suppose now that $D$ is
a pseudo-effective $\QQ$-divisor on a surface $X$.
Recall that $D$ has a Zariski decomposition
$D=P+N$, where $P$ and $N$ are $\QQ$-divisors, with $P$
nef and $N$ effective, inducing for all divisible
enough $m$ isomorphisms:
\begin{equation}\label{isom_Zariski} H^0(X,mP)\ \simeq \
H^0(X,mD)
\end{equation} (see \cite{Badescu} for details). It is
shown in
\cite{elmnp}, Example 1.11, that we have ${\rm
Supp}(N)\subseteq
\B+(D)=\B+(P)={\rm Null}(P)$.

If $C$ is an irreducible curve on $X$, and if
$C\not\subseteq{\rm Supp}(N)$, then
(\ref{isom_Zariski}) induces an equality $\vol_{X\vert
C}(D) =\vol_{X\vert C}(P)$. If, moreover, $C$ is not
contained in
${\rm Null}(P)$, then Corollary~\ref{volume_for_nef}
gives
$\vol_{X\vert C}(D)=(P\cdot C)$.
\end{example}

A theorem of Angehrn and Siu \cite{AS} on effective base-point
freeness for line bundles of the form $K_X+L$, with $L$ ample, can
be extended to the case of arbitrary big divisors, as follows.

\begin{theorem}\label{Siu}
Let $L$ be a line bundle on a smooth, projective $n$-dimensional
variety $X$. If $x\not\in\B+(L)$ is such that for every
positive-dimensional subvariety $V$ through $x$ we have
$$\vol_{X\vert V}(L)\ > \ M^{\dim(V)},$$ where
$M={{n+1}\choose 2}$, then $x$ is not in the base locus of $K_X+L$.
\end{theorem}

We do not give the proof of this statement, as it is rather
straightforward, combining the method of Angehrn and Siu (see also
Theorem~10.4.2 in \cite{positivity}) with our generalized Fujita
Approximation. We mention also that one can give similar uniform
 bounds that imply that
$K_X+L$ separates two points.

\begin{remark}\label{general_inequality} If $V$ is a
$d$-dimensional subvariety of $X$ that is contained in
$\B+(D)$ but not in $\BB(D)$, then both $\vol_{X\vert
V}(D)$ and
$\alinser{ D^d\cdot V}$ are defined, but they
are not equal in general. Example~\ref{failure} below
gives a big globally generated line bundle for which
$\vol_{X\vert C}(L)<\alinser{ L\cdot C}$ for
some curve $C$ contained in $\B+(L)$.

However, if $C$ is a curve not contained in $\BB(D)$
and $\vol_{X\vert C}(D)=0$, then we have $\parallel
D\cdot C\parallel=0$. Indeed, if $\pi_m\colon X_m\to X$
and $\pi_m^*(mD)=M_m+E_m$ are as in the definition of
asymptotic intersection numbers, then
$\widetilde{C}_m$ is not contained in the support of
$E_m$, so
$$\vol_{X\vert
C}(D)=\vol_{X_m\vert\widetilde{C}_m}(\pi_m^*(D))
\ \geq \
\frac{\vol_{X_m\vert\widetilde{C}_m}(M_m)}{m},$$ so the
right hand side is zero. Using the diagram in the proof
of Lemma~\ref{pull_back}, we deduce that
$\widetilde{C}_{m}$ is contracted by the morphism
defined by
$|M_m|$, so $(M_m\cdot \widetilde{C}_m)=0$ for every
$m$.
\end{remark}

\begin{example} Inspired by ideas of Tsuji
\cite{tsuji2}, Takayama introduced in \cite{takayama1}
the following asymptotic intersection number of a
$\QQ$-divisor $D$ with an arbitrary curve $C \subset X$:
$$\alinser{D\cdot C}' \ :=\
\limsup_{m\to\infty}\frac{{\rm deg}\Big(\OO_C(mD)\otimes
\AMI{mD}|_C\Big)}{m},$$ where the degree
of $\AMI{mD}|_C$ is defined as the
degree of the invertible sheaf
$\AMI{mD}\cdot \OO_{C'}$ on
the normalization
$\nu: C' \rightarrow C$. One can show, using the
multiplier ideal interpretation in the statement of
Theorem \ref{generalized_Fujita}, that if $C
\not\subseteq \B+(D)$, then
$\alinser{D\cdot C}'=\vol_{X\vert C}(D)$. If
$C \subseteq \B+(D)$, then we have only one inequality,
namely $\alinser{ D\cdot C}'\geq\vol_{X\vert
C}(D)$.
\end{example}

\begin{example}[Interpretation of Kodaira-Iitaka
dimension]\label{iitaka} Let $X$ be a smooth
 projective
variety of dimension $n$ and $D$ a $\QQ$-divisor on $X$ such that
its Iitaka dimension $\kappa(D)$ is nonnegative. Define $r(D)$ to be
the largest nonnegative integer $d$ such that through a very general
point on $X$ there is a $d$-dimensional irreducible subvariety $V$
with the following property: for every curve $C\subseteq V$ that is
not contained in $\BB(D)$ we have $\parallel D\cdot C\parallel=0$.
Then we have $r(D)=n-\kappa(D)$.  In order to show this, consider
$m$ divisible enough such that the rational map $\phi_m$ to $\PP^N$
defined by the complete linear series $|mD|$ satisfies
$\dim(\phi_m(X)) = \kappa(D)$. With the usual notations, we also
have a morphism $\psi_m\colon X_m \rightarrow \PP^N$ defined by the
basepoint-free linear series $|M_m|$. Let $C$ be a curve in $X$ not
contained in $\BB(D)$. If $\phi_m$ is defined at the generic point
of $C$, then $C$ is contracted by $\phi_m$ if and only if its proper
transform $\widetilde{C}_m$ is contracted by $\psi_m$. But this last
statement is equivalent to
 $(M_m\cdot \widetilde{C}_m) =0$. In particular, this
shows that if $V$ is as in the definition of
$r(D)$, then $V$ is contracted by every $\phi_m$.

Consider now the Iitaka fibration corresponding to $D$
(see \cite{positivity},
\S 2.1.C). This is a morphism of normal varieties
(defined up to birational equivalence) $\phi_{\infty}
\colon X_{\infty}\to Y_{\infty}$ having connected
fibers and such that for all $m$ divisible enough we
have a commutative diagram
\[
\begin{CD} X_{\infty}@>{u_{\infty}}>> X \\
@V{\phi_{\infty}}VV @VV{\phi_m}V \\ Y_{\infty}
@>u_m>>Y_m
\end{CD}
\] with $u_{\infty}$ a birational morphism and $u_m$ a
birational map. Let $U$ be the set of points $x$ in
$X\smallsetminus\BB(D)$ such that $u_{\infty}$ is an isomorphism
over a neighborhood of $x$ and $\phi_{\infty}(u_{\infty}^{-1}(x))$
lies in an open subset on which $u_m$ is an isomorphism. We see that
if $C$ is a curve in $X$ that intersects $U$ and
$\widetilde{C}\subseteq X_{\infty}$ is its proper transform, then
$C$ is contracted by some $\phi_m$ (with $m$ divisible enough) if
and only if $\widetilde{C}$ is contracted by $\phi_{\infty}$. Since
this condition is independent on $m$, we see that it is satisfied if
and only if $\alinser{D\cdot C}=0$. Our formula for $r(D)$ now
follows easily.
\end{example}

We restate the conclusion of the above discussion as
follows. Compare this with Theorem~1.3 in
\cite{takayama1}, where a similar result is proved with
$\alinser{D\cdot C}$ replaced by
$\alinser{ D\cdot C}'$.

\begin{corollary}\label{iitaka_equivalence} If
$\phi_{\infty}\colon X_{\infty}\to Y_{\infty}$ is the
Iitaka fibration corresponding to $D$ and if $C$ is a
curve through a very general point on $X$, then its
proper transform $\widetilde{C}
\subseteq X_{\infty}$ is contracted by $\phi_{\infty}$
if and only if
$\alinser{ D\cdot C}=0$. The analogous
assertion holds if we replace the morphism
$\phi_{\infty}$ by the rational map on $X$ defined by
$|mD|$, with $m$ divisible enough.
\end{corollary}

\section{The proof of Generalized Fujita Approximation}

For the proof of Theorem \ref{generalized_Fujita}, we
will need a few lemmas.

\begin{lemma}\label{bounded_difference} For any big
$\QQ$-divisor D, if $\frb_m$ denotes the ideal defining
the base locus of $|mD|$, then there exists an
effective integral divisor $G$ on $X$ -- which one may
take to be very ample --  such that
$$\frb_m(-G) \ \subseteq \ \AMI{mD}(-G)
\ \subseteq \ \frb_m,$$ for all $m$ sufficiently large
and divisible. If $V\not\subseteq \B+ (D)$, then
$G$ can be chosen such that $V\not \subseteq {\rm
Supp}(G)$.
\end{lemma}
\begin{proof} The first statement appears in
\cite{positivity} Theorem 11.2.21, but we recall the
construction in order to emphasize the second
point. Specifically, let
$H$ be a very ample line bundle on $X$, and consider
the ample line bundle $A:= K_X + (n+1)H$, where $n$ is
the dimension of $X$. For $a\gg 0$ sufficiently
divisible, $\OO_X(aD - A)$ is a big line bundle with sections, and
 for any $G$ in
$|aD-A|$ the sequence of inclusions in the Lemma holds.

Now note that since $a$ is large and
divisible enough,
${\rm Bs}(|aD - A|)$ is contained in $\B+(D)$. Indeed,
if $m$ is large enough then
$\B+(D)=\BB(D-\frac{1}{m}D)$, and if $a$ is divisible
enough, then
$$\BB(D-\frac{1}{m}A)\ =\ {\rm
Bs}(|aD-\frac{a}{m}A|)_{\rm red}\ \supseteq \ {\rm
Bs}(|aD-A|)_{\rm red}.$$ Thus if
$V\not\subseteq \B+ (D)$, then there exists $G\in
|aD-A|$ such that $V\not \subseteq {\rm Supp}(G)$.
\end{proof}

\begin{lemma}\label{semiample} Let $L$ be a big
semiample line bundle on $X$ and suppose that $V$ is a
$d$-dimensional subvariety of $X$ such that
$V\not\subseteq \B+ (L)$. Then
$\vol_{X\vert V}(L) = (L^d\cdot V).$
\end{lemma}
\begin{proof} By choosing $k\gg 0$, we may assume that
the morphism $f:X \rightarrow \PP^N$ given by the
linear series $|kL|$ is birational onto its image, with
trivial Stein factorization. Since $V\not\subseteq \B+
(L)$, we may also assume that the restriction of $f$ to
V is birational onto its image. The argument in the
proof of Lemma~\ref{pull_back} reduces us to the case
of a very ample line bundle, when the equality is clear.
\end{proof}

\begin{lemma}\label{scale_multiplier} If $V$ is a
subvariety of $X$ of dimension $d\geq 1$, and if $D$ is
an integral divisor on $X$, then for every $k$ we have
\begin{equation}\label{same_limit}
\limsup_{m\to\infty}\frac{\hh{0}{V}{\OO(mD)\otimes
\AMI{mD}|_V}}{m^d} \ = \
\limsup_{p\to\infty}\frac{\hh{0}{V}{\OO(pkD)\otimes
\AMI{pkD}|_V}}{k^dp^d} .
\end{equation}
\end{lemma}

\begin{proof} Let's denote the left hand side in
(\ref{same_limit}) by $L_1$ and the right hand side by
$L_2$. We obviously have $L_1 \geq L_2$ and we need to
prove the reverse inequality. To this end let $A$ be a
very ample line bundle on $V$, and fix a very general
divisor $H \in \linser{A}$. Assuming as we may that $H$
doesn't contain any of the subvarieties defined by
the associated primes of   the ideal sheaves
$\AMI{mD}|_V$, we have for every $m \ge 0$ an exact
sequence
\begin{multline*}
0 \lra \OO_V(mD) \otimes \AMI{mD}|_V
\overset{\cdot H}{\lra}
 \OO_V(mD + A) \otimes
\AMI{mD}|_V \\  \lra   \OO_H(mD + A) \otimes \AMI{mD}|_H
\lra 0. \end{multline*}
Since in any event
\[ \hh{0}{H} {\OO_H(mD + A) \otimes \AMI{mD}|_H} \ \le
\ \hh{0}{H}{\OO_H(mD + A)} \ = \ O(m^{d-1}), \]
we see in the first place that
\begin{equation}  L_1 \ = \
\limsup_{m\to\infty}\frac{\hh{0}{V}{\OO(mD+A)\otimes
\AMI{mD}|_V}}{m^d}. \tag{*}
\end{equation}

Now given $k \ge 1$ and $m \gg 0$, write $m = pk +
\ell$ with $0 \le \ell \le k-1$. Choose a very ample
line bundle
$A$ which is sufficiently positive so that $A + \ell D$
is very ample for each $0 \le \ell \le k -1$, and fix a
very general divisor $H_\ell \in \linser{\ell D + A}$.
As above we have an exact sequence
\begin{multline*}
0 \lra \OO_V(pkD) \otimes \AMI{pkD}|_V \overset{\cdot
H_\ell}{\lra} \OO_V(mD +A ) \otimes \AMI{pkD}|_V \\ \lra
\OO_{H_\ell}(mD +A ) \otimes \AMI{pkD}|_{H_\ell}\lra 0.
\end{multline*}
Noting that $\AMI{mD} \subseteq \AMI{pkD}$, we find
as before that
\begin{align*}
L_2 \
&= \ \limsup_{p\to\infty}\frac{\hh{0}{V}{\OO(pkD)\otimes
\AMI{pkD}|_V}}{p^dk^d} \\ &= \
\limsup_{m\to\infty}\frac{\hh{0}{V}{\OO(mD+A)\otimes
\AMI{(m-\ell)D}|_V}}{m^d} \\ & \ge \
\limsup_{m\to\infty}
\frac{\hh{0}{V}{\OO(mD+A)\otimes
\AMI{mD}|_V}}{m^d}.
\end{align*}
Together with (*), this gives the required inequality
$L_2 \ge L_1$.
\end{proof}

\begin{proof}[Proof of Theorem~\ref{generalized_Fujita}]
Consider a common log resolution $\pi_m: X_m \ra X$ for the ideal
$\frb_m$ defining the base locus of $|mD|$ and for $\ci(\parallel
mD\parallel)$. We denote $\frb_m\cdot \OO_{X_m} = \OO _{X_m}(-E_m)$
and $\ci(\parallel mD\parallel)\cdot \OO_{X_m} = \OO _{X_m}(-F_m)$,
and also $M_m := \pi_m^*(mD) - E_m$ and $N_m: =\pi_m^*(mD) - F_m$.
Since $\frb_m \subseteq \ci(\parallel mD\parallel)$, we have $E_m
\geq F_m$. Note also that $|M_m|$ is the moving part of the linear
series $|\pi_m^* (mD)|$, so it  is basepoint-free. Moreover, since
$V\not\subseteq\B+(D)$, we have $\widetilde{V}_m\not
\subseteq\B+(M_m)$.

We have, to begin with
$$\vol _{X
\vert V}(D) \geq
\frac{\vol_{X_m\vert\widetilde{V}_m}(M_m)}{m^d}=
\frac{\vol(M_m|_{\widetilde{V}_m})}{m^d} =  \frac{
M_m^d\cdot \widetilde{V}_m}{m^d},$$ where the first
inequality follows easily from the definition and
Lemma~\ref{pull_back}, while the second and third
equalities follow from Lemma \ref{semiample}. This
implies that
$$\vol_{X\vert V}(D) \geq ~\alinser{ D^d \cdot
V}.$$ On the other hand, since for any
(divisible enough) $m$ we have
$$\HH{0}{X}{\OO(mD)} \ = \  \HH{0}{X}{\OO(mD)\otimes
\AMI{mD}},$$ we see immediately that
$$\limsup_{m\to\infty}\frac{\hh{0}{V}{ \OO(mD)\otimes
\ci(\parallel mD\parallel)|_V)}}{m^d/d!}
\geq\vol_{X\vert V}(D).$$ To finish, it suffices then
to prove
$$\alinser{ D^d \cdot V} \geq
\limsup_{m\to\infty}\frac{h^0(V, \OO(mD)\otimes
\AMI{mD}|_V)}{m^d/d!}.$$

To this end, we first apply the inclusions
$$\frb_m(-G)\subseteq \AMI{mD}(-G)\subseteq
\frb_m$$ given by Lemma
\ref{bounded_difference}. On $X_m$, this immediately
implies:
\begin{equation}\label{comparing_volumes}
\vol ({(M_m + \pi_m^*G)}|_{\widetilde{V}_m})\geq \vol
({N_m}|_{\widetilde{V}_m})
\geq \vol ({M_m}|_{\widetilde{V}_m}).
\end{equation}

On the other hand, by pulling back to $\widetilde{V}_m$
in two different ways, we have:
$$h^0 (\widetilde{V}_m, kN_m|_{\widetilde{V}_m}) \geq h^0(V,
\OO(kmD)\otimes
\AMI{mD}^k|_V) \geq h^0(V,
\OO(kmD)\otimes
\AMI{ kmD}|_V),$$ where for the last
inequality we use the Subadditivity Theorem (see
\cite{del} or \cite{positivity}, Theorem 11.2.3). Thus
by multiplying the inequalities by $d!/(km)^d$ and
letting first $k$ and then $m$ go to $\infty$, we
obtain:
\begin{equation}\label{inequality_1}
\limsup_{m\to\infty}\frac{\vol
({N_m}|_{\widetilde{V}_m})}{m^d} \geq
\limsup_{m\to\infty}\frac{h^0(V, \OO(mD)\otimes
\AMI{ mD}|_V)}{m^d/d!},
\end{equation} thanks to Lemma \ref{scale_multiplier}.

Combining (\ref{comparing_volumes}) and
(\ref{inequality_1}), we are then done if we show
\begin{equation}\label{formula11}
\limsup_{m\to\infty}\frac{\vol ((M_m + \pi_m^*
G)|_{\widetilde{V}_m})}{m^d} =
\limsup_{m\to\infty
}\frac{\vol
({M_m}|_{\widetilde{V}_m})}{m^d}.
\end{equation}
But since the bundles in question are nef, the volumes
appearing here are computed as intersection numbers.
Expanding out the one on the left, we find that
\eqref{formula11} will follow if we show that
\begin{equation}\label{Ineq.of.Int.Nos}
\limsup_{m \to \infty} \frac{ \big( M_m^{d-i} \cdot
(\pi_m^* G)^i \cdot {\widetilde{V}_m} \big)}{m^d} \ = \
0
\end{equation}
when $i > 0$.
But we can find a fixed ample line bundle $A$ on $X$
such that $\pi^*(mA) - M_m$ is effective (simply take
$A$ such that $A - L$ is effective). Then
\[
\big( M_m^{d-i} \cdot
(\pi_m^* G)^i \cdot {\widetilde{V}_m} \big)\ \le \
\big( (\pi_m^*(mA)^{d-i} \cdot (\pi_m^* G)^i \cdot
{\widetilde{V}_m} \big),
\]
and \eqref{Ineq.of.Int.Nos} follows.
\end{proof}

\section{Approximating restricted volumes via jet
separation}  We start with a simple lemma which shows that
separating jets at a set of points on a subvariety gives a lower
bound for the restricted volume. If $x$ is a point on a subvariety
$W$ of $X$, we denote by $\frmm_{W,x}$ the ideal defining $x$ in the
local ring $\OO_{W,x}$ of $W$ at $x$. If $x_1,\ldots,x_N$ are points
on $X$, then we say that a line bundle $L$ on $X$ simultaneously
separates $p_i$-jets at each $x_i$ if the map
$$\HH{0}{X}{L}\longrightarrow\bigoplus_{i=1}^N
\HH{0}{X}{L\otimes\OO_{X,x_i}/ \frmm_{X,x_i}^{p_i+1}}$$ is
surjective.

\begin{lemma}\label{volume_from_jets}
 Let $L$ be a line bundle on the smooth variety $X$,
and let $V\subseteq X$ be a subvariety of dimension $d\geq 1$. If
$x_1,\ldots,x_N$ are distinct points on $V$ such that $L$
simultaneously separates $p_i$-jets on $X$ at each of the points
$x_i$, then
$$\vol_{X\vert V}(L)\geq\sum_{i=1}^N{\rm
mult}_{x_i}V\cdot p_i^d.$$
\end{lemma}

\begin{proof} If $L$ separates $p_i$-jets at each
$x_i$, then by taking polynomials in sections one sees that if $m
\ge 1$ then $mL$ separates $mp_i$-jets on $X$ at $x_i$ for all $i$.
Moreover, for every such $m$ we have a commutative diagram
\[
\begin{CD} H^0(X,mL)
@>>>\bigoplus_iH^0(X,mL\otimes\OO_{X,x_i}/\frmm_{X,x_i}^{mp_i+1})
\\ @V
VV @VVV \\ H^0(V,mL\vert_V)
@>>>\bigoplus_iH^0(V,mL\otimes\OO_{V,x_i}/
\frmm_{V,x_i}^{mp_i+1}).
\end{CD}
\]
As the right vertical map is evidently surjective, we
deduce that
$$\dim_{\CC}\ {\rm Im}\HH{0}{X|V}{mL}\
\geq \
 \sum_{i=1}^N\dim_{\CC}\OO_{V,x_i}/
\frmm_{V,x_i}^{mp_i+1},$$ and the lemma follows.
\end{proof}

We now prove the converse to Lemma \ref{volume_from_jets}, showing
that we can approximate restricted volumes by separation of jets at
general points on the subvariety. It is convenient to make the
following definition. Let $D$ be a $\QQ$-divisor on the smooth
projective variety $X$, and let $V$ be a subvariety of $X$ of
dimension $d\geq 1$. For every positive integer $N$, let
$\epsilon_V(D,N)$ be the supremum of the set of nonnegative rational
numbers $t$ with the property that for some $m$ with $mt\in\ZZ$ and
$mD$ an integral divisor, $mD$ simultaneously separates $mt$-jets at
every general set of points $x_1,\ldots,x_N\in V$ (if there is no
such $t$, then we put $\epsilon_V(D,N)=0$). Note that with this
notation, Lemma~\ref{volume_from_jets} implies $\vol_{X\vert
V}(D)\geq N\cdot\epsilon_V(D,N)^d$.

\begin{theorem}\label{separation_of_jets} If $D$ is a
$\QQ$-divisor on the smooth projective variety $X$, and $V\subseteq
X$ is a subvariety of dimension $d\geq 1$ such that $V\not\subseteq
\B+(D)$, then
$$\sup_{N\geq 1}N\cdot\epsilon_V(D,N)^d=\limsup_{N\to\infty}
N\cdot\epsilon_V(D,N)^d=\vol_{X\vert V}(D).$$
\end{theorem}

\begin{remark} The above statement is
inspired by a key step in the proof of Fujita's Approximation
Theorem from \cite{nakamaye2}. In \emph{loc. cit.} one proves a
variant of this statement when $V=X$. However, we take the opposite
approach, and we deduce Theorem~\ref{separation_of_jets} from our
generalization of Fujita's Theorem to restricted volumes.
\end{remark}

\begin{proof}[Proof of Theorem~\ref{separation_of_jets}]
We need to show that for every $\delta>0$ we can find arbitrarily
large values of $N$ such that for a suitable positive
$\epsilon\in\QQ$ with $N\epsilon^d>\vol_{X\vert V}(D)-\delta$, and
for some positive integer $m$ such that $m\epsilon\in\ZZ$ and $mD$
is an integral divisor, $mD$ simultaneously separates
$m\epsilon$-jets on $X$ at any general points $x_1,\ldots,x_N\in V$.
Suppose first that we know this when $D$ is ample. Since
$V\not\subseteq\B+(D)$, it follows from
Theorem~\ref{generalized_Fujita} and
Proposition~\ref{general_resolution} that we can find a proper
morphism $\pi : X'\longrightarrow X$ that is an isomorphism over the
generic point of $V$, and a decomposition $\pi^*D=A+E$, with $A$
ample, $E$ effective and $\widetilde{V}\not\subseteq{\rm Supp}(E)$
such that $(A^d\cdot \widetilde{V})>\vol_{X\vert V}(D)-\delta/2$ (we
have denoted by
 $\widetilde{V}$ the proper transform of $V$). We apply
the ample case for $A$, $\widetilde{V}$ and $\delta/2$ to get
$\epsilon$, $N$ and $m$. We may clearly assume also that $mE$ is
integral. If $x_1,\ldots,x_N\in \widetilde{V}$ are general points
(in particular they do not lie on the union of the support of $E$
with the exceptional locus of $\pi$), and if we identify the $x_i$
with their projections to $X$, then we have a commutative diagram

\[
\begin{CD} H^0(X',\OO(mA)) @>>\psi_1> \bigoplus_{i=1}^N
H^0(X',\OO(mA)\otimes\OO_{X',x_i}/
\frmm_{X',x_i}^{m\epsilon+1})\\ @VV\phi_1V @VV\phi_2V\\
H^0(X,\OO(mD)) @>>\psi_2> \bigoplus_{i=1}^N
H^0(X,\OO(mD)\otimes\OO_{X,x_i}/
\frmm_{X,x_i}^{m\epsilon+1})\\
\end{CD}
\] where $\phi_1$ and $\phi_2$ are induced by
multiplication with the section defining $mE$. Hence $\phi_2$ is an
isomorphism, and since $\psi_1$ is surjective, $\psi_2$ is
surjective, too. Therefore we get our statement for $D$, $V$ and
$\delta$.

It follows that in order to prove the theorem we may
assume that $D$ is ample, in which case
$\vol_{X\vert V}(D)=(D^d\cdot V)$. Moreover, by
replacing $D$ with a suitable power, we may suppose
that it is very ample.

We make a parenthesis to recall the following well-known fact.
Suppose that $L$ is an ample line bundle on a variety $X$, and
suppose that $\Gamma=\{x_1,\ldots,x_M\}$ is a set of smooth points
on $X$. Let $f: X'\longrightarrow X$ be the blowing-up along
$\Gamma$ with exceptional divisor $F=\sum_{i=1}^MF_i$. If $\beta>0$,
then $f^*(L)-\beta F$ is nef if and only if for every positive
rational number $\epsilon<\beta$, if $k$ is divisible enough, then
$kL$ separates $k\epsilon$-jets at $x_1,\ldots,x_M$. Moreover, this
is the case if and only if for every irreducible curve $C$ in $X$ we
have
$$(L\cdot C)\ \geq\ \beta\sum_{i=1}^
M{\rm
mult}_{x_i}(C).$$ If this holds, and $\beta'<\beta$,
then $f^*(L)-\beta' F$ is ample on $X'$. Note also that
if the condition is satisfied for $\Gamma$, then it is
satisfied for any subset $\Gamma'\subseteq\Gamma$ too.
We take $p$ such that $V$ is cut out by equations in
$|pD|$. Let
$H_{d+1},\ldots,H_n\in |pD|$ be general elements
vanishing on $V$, so $V$ is an irreducible component of
$W:=\bigcap_{i=d+1}^nH_i$. Moreover,
$V=W$ scheme-theoretically at the generic point of $V$,
and $W\setminus V$ is smooth of dimension $d$. Let
$H_1,\ldots,H_d$ be general elements in $|pD|$, so the
following sets
$$\Gamma'\ :=\ V\cap H_1\cap\ldots\cap H_d\ \subseteq \
\Gamma:=\bigcap_{i=1}^nH_i$$ are smooth and
zero-dimensional. Let $x_1,\ldots,x_M$ be the points in
$\Gamma$, and suppose that they are numbered such that
the first $N$ are the points in $\Gamma'$, where
$N=p^d(D^d\cdot V)$.  Note that if $C\subseteq X$ is an
irreducible curve, then there is $i\leq n$ such that
$C\not\subseteq H_i$. It follows from B\'{e}zout's
theorem that
$$(H_i\cdot C)\ \geq\ \sum_{j=1}^M{\rm mult}_{x_j}(C).$$
We deduce from the previous discussion that if $\epsilon<1/p$ and
$m$ is divisible enough, then $mD$ separates $m\epsilon$-jets at
$x_1,\ldots,x_M$, hence at $x_1,\ldots,x_N$.  Given $\delta$ and
$\eta>0$, we choose $p\gg 0$ as above and such that $1/p<\eta$. If
$\epsilon$ is such that $1>(p\epsilon)^d>1- \frac{\delta}{(D^d\cdot
V)}$, we see that for $m$ divisible enough the points
$x_1,\ldots,x_N$ satisfy our requirement. It is now standard (using
the behavior of ampleness in families) to deduce that the same
property holds for any general set of points in $V$. We end by
noting that since $p$ can be taken arbitrarily large, the same is
true for $N$.
\end{proof}

For future reference, we recall the following well-known facts.

\begin{remark}\label{extra_divisor}
Suppose that $B$ is an ample $\QQ$-divisor and $x_1,\ldots,x_N\in X$
are such that for $m$ divisible enough, $mB$ simultaneously
separates $m\epsilon'$-jets at $x_1,\ldots,x_N$. If
$\epsilon<\epsilon'$ and if $m$ is divisible enough, then the linear
system
$$\{P\in |mB|\mid
 {\rm ord}_{x_i}(P)\geq m\epsilon+1\,{\rm for}\,{\rm all}\,i\}$$ induces
a base-point free linear system on
$X\smallsetminus\{x_1,\ldots,x_N\}$.

Indeed, if $\pi\colon X'\to X$ is the blowing-up at $x_1,\ldots,x_N$
with exceptional divisors $E_1,\ldots,E_N$, our hypothesis implies
that $B':=\pi^*B-\epsilon(E_1+\ldots+E_N)$ is ample. Hence every $m$
such that $mB'$ is integral and globally generated satisfies our
requirement.
\end{remark}

\begin{remark}\label{separation2}
Let $B, x_1,\ldots,x_N$ and $\epsilon$, $\epsilon'$ be as in the
previous remark. For every integral divisor $H$ on $X$, if $m$ is
divisible enough, then $(mB-H)$ separates $m\epsilon$-jets at
$x_1,\ldots,x_N$. Indeed, if $\pi$ is as before, then
$\pi^*B-\epsilon(E_1+\ldots+E_N)$ is ample. It follows that if $m$
is divisible enough, then $mB-H$ is an integral divisor,
$m\epsilon\in\ZZ$ and $\pi^*(mB-H)-m\epsilon(E_1+\ldots+E_N)$ is
globally generated. This implies that $(mB-H)$ simultaneously
separates $m\epsilon$-jets at $x_1,\ldots,x_N$.
\end{remark}

\begin{remark}\label{separation3}
Suppose now that $D$ and $V$ are as in
Theorem~\ref{separation_of_jets}, and that $H$ is an integral
effective divisor on $X$. If $N\geq 1$ and
$\epsilon<\epsilon_V(D,N)$, then for $m$ divisible enough both $mD$
and $mD-H$ simultaneously separate $m\epsilon$-jets on $X$ at every
general set of points $x_1,\ldots,x_N$ in $V$. Indeed, we argue as
in the proof of the theorem: we use Theorem~\ref{generalized_Fujita}
to reduce ourselves to the case of an ample line divisor $A$ on some
model $X'$ over $X$. We apply for $A$ the argument in the previous
remark, and use the fact that if $(mA-\pi^*(H))$ separates jets,
then so does $mD-H$.
\end{remark}

\section{Components of $\BB_+$ and the restricted
volume function}

Given an $\RR$-divisor $D$, we have $\B+(D')\subseteq\B+(D)$ for
every $\RR$-divisor $D'$ in a suitable open neighborhood of $D$. It
follows that given a subvariety $V$ of $X$, the set $\{D\mid
V\not\subseteq\B+(D)\}$ is an open subset of the big cone.

\begin{definition}
Given a subvariety $V\subseteq X$, we denote by ${\rm
Big}^V(X)_{\RR}$ the set of big $\RR$-divisor classes $D$ such that
$V$ is not a proper subset of an irreducible component of $\B+(D)$.
It is clear that ${\rm Big}^V(X)_{\RR}$ is an open convex subcone of
the big cone.
\end{definition}

The behavior of the restricted volumes as functions on subsets of
the big cone can be summarized in the following theorem.

\begin{theorem}\label{continuity} (a) If $V$ is a fixed
subvariety of a smooth projective variety $X$, with
$d=\dim(V)\geq 1$, then the map
$$\xi\longrightarrow\vol_{X\vert V}(\xi)$$ defined on
the set of $\QQ$-divisor classes $\xi$ such that $V\not\subseteq
\B+(\xi)$ is continuous and can be extended to a continuous function
on the open set of all such $\RR$-divisor classes. Moreover, it
satisfies the concavity relation
$$\vol_{X\vert V}(D_1+D_2)^{1/d}\geq\vol_{X\vert
V}(D_1)^{1/d} +\vol_{X\vert V}(D_2)^{1/d}$$ for every
$\RR$-classes as above.

(b) If $D$ is a $\QQ$-divisor such that $V$ is an
irreducible component of
$\B+(D)$, then $\vol_{X\vert V}(D)=0$. Moreover, if we
put $\vol_{X\vert V}(\xi)=0$ for every $\xi\in {\rm
Big}^V(X)_{\RR}$ such that $V\subseteq \B+(\xi)$, then
the function
$\xi\longrightarrow\vol_{X\vert V}(\xi)$ is continuous
over the entire cone ${\rm Big}^V(X)_{\RR}$.
\end{theorem}

\begin{proof} The proof of part (a) is quite standard,
and we present it in what follows. Part (b) is the main
technical result of the paper, and we present its proof
separately (cf. Theorem \ref{main} below).

Fix
 ample $\QQ$-divisors $A_1,\ldots,A_r$ whose classes
in
$N^1(X)_{\QQ}$ form a basis. It is convenient to take
on $N^1(X)_{\QQ}$ the norm
$\parallel\sum_{i=1}^r\alpha_iA_i\parallel
=\max_i|\alpha_i|$.  For a positive rational number
$s$, and for a $\QQ$-divisor $D_0$ such that
$V\not\subseteq\B+(D_0-s\sum_{i=1}^rA_i)$, we denote by
$T(D_0,s)$ the set of divisor classes
$D_0-(p_1A_1+\ldots p_rA_r)$, where $0\leq p_i\leq s$
are rational numbers. It is enough to show that for
every such box $T(D_0,s)$ there is a constant
$C=C(D_0,s)$ such that
\begin{equation}\label{suf1} |\vol_{X\vert
V}(D_1)-\vol_{X\vert V}(D_2)|
\leq C\cdot\parallel D_1-D_2\parallel,
\end{equation} for every $D_1$ and $D_2$ in $T(D_0,s)$.

Let $H$ be a fixed ample $\QQ$-divisor, and $D$ in
$T(D_0,s)$. If $\epsilon_0$ is such that
$V\not\subseteq\B+(D_0-s\sum_iA_i-\epsilon_0H)$ and if
$\epsilon\leq\epsilon_0$, then
\begin{equation}\label{ineq1}
\vol_{X\vert V}(D-\epsilon H) \ge \vol_{X\vert
V}((1-\epsilon/\epsilon_0)D) =
(1-\epsilon/\epsilon_0)^d \vol_{X\vert V}(D).
\end{equation}

By the openness of the ample cone, there is a positive
real number $b$, such that if $E$ is a $\QQ$-divisor
with $\parallel E \parallel \le b$, then $H- E$ is
ample. If $A$ is a $\QQ$-divisor such that
$\parallel A \parallel \le b\cdot \epsilon_0$, we have
that
$\frac{\parallel A \parallel}{b}H-A$ is ample. Combined
with (\ref{ineq1}), this shows that for every $D\in
T(D_0,s)$
\begin{equation}\label{ineq2}
\vol_{X\vert V}(D - A) \geq \vol_{X\vert V}
\left(D-  \frac{\parallel A \parallel}{b} H\right)
\geq \left(1-  \frac{\parallel A
\parallel}{b\epsilon_0}\right)^d
\vol_{X\vert V}(D).
\end{equation} As $d\geq 1$, we see that there is a
constant $C'$ such that for every $D$ and $A$ as above
we have
\begin{equation}\label{ineq111}
\vol_{X\vert V}(D)- \vol_{X\vert V}(D - A)
\le  C'\cdot \parallel A\parallel \cdot \vol_{X\vert
V}(D)\leq C'\cdot
\vol_{X\vert V}(D_0)\cdot\parallel A\parallel.
\end{equation}

Suppose now that $D\in T(D_0,s)$ and that $A$ is an effective linear
combination of the $A_i$ such that $D-A\in T(D_0,s)$. If $m\gg 0$
and if we apply (\ref{ineq111}) successively to $D-\frac{i}{m}A$ and
$\frac{1}{m}A$ for $0\leq i\leq m-1$, we deduce
\begin{equation}\label{ineq112} |\vol_{X\vert
V}(D)-\vol_{X\vert V}(D-A)|
\leq \frac{C}{2}\cdot\parallel A\parallel,
\end{equation} where $C=2C'\cdot\vol_X(V,D_0)$.

We finish the proof as in the case of the usual volume
function (see \cite{positivity}, Theorem 2.2.44). Note
that if $D_1$, $D_2\in T(D_0,s)$, then we may write
$D_2=D_1+E-F$, where $E$ and $F$ are effective linear
combinations of the $A_i$. We apply (\ref{ineq112}) to
get
$$\vol_{X\vert V}(D_1)-\vol_{X\vert V}(D_1-F)
\leq \frac{C}{2}\cdot\parallel F\parallel,$$
$$\vol_{X\vert V}(D_2)-\vol_{X\vert V}(D_1-F)
\leq \frac{C}{2}\cdot\parallel E\parallel.$$ Since
$\parallel D_1-D_2\parallel =\max\{\parallel
E\parallel, \parallel F\parallel\}$, (\ref{suf1})
follows from the triangle inequality.

The fact that $\vol_{X\vert V}(-)^{1/d}$ is concave on
the set of $\RR$-divisor classes $D$ such that
$V\not\subseteq\B+(D)$ follows by continuity from
Corollary~\ref{concavity}.
\end{proof}

\begin{remark} Note that in (\ref{ineq111}) in the
above proof, we may take $A$ to be numerically trivial.
This gives another way of seeing that if
$V\not\subseteq\B+(D)$, then the restricted volume
$\vol_{X\vert V}(D)$ depends only on the numerical
class of $D$.
\end{remark}

\begin{remark}
We chose to give the above proof of Theorem~\ref{continuity} a) that
shows that more generally, every homogeneous function defined on the
rational points of an open convex cone containing the ample cone,
and which is non-decreasing with respect to adding an ample class is
locally Lipschitz continuous. Alternatively, the assertion in the
theorem could be deduced from the more subtle concavity property of
the restricted volume function, plus the following well-known fact:
a homogeneous convex function defined on the rational points of a
convex domain is Lipschitz around every point in the domain (in
particular, it is locally uniformly continuous, and therefore it can
be extended by continuity to the whole domain).
\end{remark}

\begin{example}\label{big_and_nef} It follows from
Theorem~\ref{continuity} that if $D$ is a nef
$\RR$-divisor and if $V\not\subseteq\B+(D)$ is a
subvariety of dimension $d\geq 1$, then
\begin{equation}\label{volume_for_nef_R}
\vol_{X\vert V}(D)=(D^d\cdot V).
\end{equation} Indeed, if $D$ is a $\QQ$-divisor, then
the assertion follows from Corollary \ref{volume_for_nef}, and the
general case follows by continuity. Moreover, the second part of
Theorem~\ref{continuity} and the continuity of the intersection form
imply that the equality (\ref{volume_for_nef_R}) still holds if
$V\subseteq\B+(D)$ is an irreducible component, and in this case
both numbers are zero (recall that each irreducible component of
$\B+(D)$ has positive dimension by Proposition~\ref{no_points}).

On the other hand, since $D$ is nef, if $V\subseteq X$
is a subvariety of dimension $d\geq 1$, we have
$D\vert_V$ big if and only if $(D^d\cdot V) >0$. This
gives the following description of $\B+(D)$, which
generalizes to the case of an $\RR$-divisor the main
result of \cite{nakamaye1}.
\end{example}

\begin{corollary}\label{cor_nakamaye} If $D$ is a nef
$\RR$-divisor, then
$\B+(D)$ is the union ${\rm Null}(D)$ of those
subvarieties $V$ of $X$ such that the restriction
$D\vert_V$ is not big.
\end{corollary}

As mentioned above, the main result on the behavior of
the restricted volume is part (b) of Theorem
\ref{continuity}. We discuss it in what follows.
\begin{theorem}\label{main} Let $X$ be a smooth
projective complex variety, and let $D$ be an
$\RR$-divisor on $X$. If $V$ is an irreducible component
of $\B+(D)$, then
\begin{equation}\label{limit}
\lim_{D'\to D}\vol_{X\vert V}(D')=0,
\end{equation} where the limit is over $\QQ$-divisors
$D'$ whose classes go to the class of $D$.
\end{theorem}

\begin{remark} Using the definition of $\vol_{X\vert
V}(-)$ on ${\rm Big}^V(X)_{\RR}$ and the concavity of
restricted volumes in Corollary~\ref{concavity}, we see
that if $V$ is an irreducible component of $\B+(D)$,
then Theorem~\ref{main} gives
$$\lim_{\xi\to D}\vol_{X\vert V}(\xi)=0,$$ where the
limit is over all $\xi\in{\rm Big}^V(X)_{\RR}$ such that
$\xi$ goes to the class of $D$. Hence we get the
assertion in part (b) of Theorem~\ref{continuity}.
\end{remark}

\begin{corollary}
\label{numerical_characterization} For any
$\QQ$-divisor $D$, the irreducible components of
$\B+(D)$ are precisely the maximal (with respect to
inclusion) $V\subseteq X$ such that $\vol_{X\vert V}(D)
= 0$.
\end{corollary}

\begin{proof} If $V$ is a component of $\B+(D)$,
Theorem \ref{main} implies that $\vol_{X\vert V}(D)=
0$. On the other hand, if $V\not\subseteq \B+ (D)$,
then we clearly
have $\vol_{X\vert V}(D) > 0$.
\end{proof}

\begin{example}[Failure of Theorem \ref{main} for
non-components]\label{failure} We give an example of a subvariety
$V$ that is properly contained in an irreducible component of the
augmented base locus and for which the conclusion of Theorem
\ref{main} is no longer true. Consider $\pi: X \rightarrow \PP^3$ to
be the blow-up of $\PP^3$ along a line $l$, with exceptional divisor
$E\cong \PP^1 \times \PP^1$ (where $\pi$ induces the projection onto
the first component). The line bundle $L:= \pi^* \OO_{\PP^3}(1)$ is
big and globally generated
 on $X$, and so by Corollary~\ref{cor_nakamaye} we have
$\B+ (L) = {\rm Null}(L) =E$.

Consider now a smooth curve $C$ of type $(2,1)$ in $E$.
It is easy to see that for all $m$, the image of the map
$$H^0 (X, mL) \longrightarrow H^0(C, mL\vert_C)$$ is
isomorphic to the image of
$$H^0(\PP^3, \OO_{\PP^3}(m)) \longrightarrow H^0(l,
\OO_{l}(m)),$$ so $\vol_{X\vert C}(L) = 1$. In
particular, it is nonzero. We also see that in this case
$\parallel L\cdot C\parallel\neq\vol_{X\vert C}(L)$
(compare with Theorem~\ref{generalized_Fujita}).
Indeed, since $L$ is globally generated we have
$\parallel L\cdot C\parallel=(L\cdot C)=2$.

Note that $\vol_{X\vert C}(-)$ is not continuous at $L$.
Indeed, let $L_m: = L- \frac{1}{m} E$ be  a sequence of
$\QQ$-divisors converging to $L$. We see that $L_m$ is
ample for $m$ large enough, which implies that
$\vol_{X\vert C}(L_m) = (L_m \cdot C)$ by Example
\ref{ample}. An easy computation shows
$$\vol_{X\vert C}(L_m) = 2 + \frac{1}{m}
\longrightarrow 2 \neq 1 = \vol_{X\vert C}(L).$$
\end{example}

\begin{remark}[Reduction to the case of positive
perturbations]\label{reduction} Note that if $D$ and
$A$ are $\QQ$-divisors and if $A$ is ample, then
$\vol_{X\vert V}(D)\leq\vol_{X\vert V}(D+A)$. It
follows that in order to prove Theorem~\ref{main}, it
is enough to consider the limit over those
$\QQ$-divisor classes $D'$ such that $D'-D$ is ample.
Indeed, suppose that $A_1,\ldots,A_r$ are ample
$\QQ$-divisors whose
 classes give a basis of
$N^1(X)_{\RR}$. We consider on $N^1(X)_{\RR}$ the norm
given by
$$\parallel\sum_i\alpha_iA_i\parallel:=\max_i|\alpha_i|.$$
It follows that given a
$\QQ$-divisor $D'\not\equiv D$, there exists a
$\QQ$-divisor $D''$ such that $\parallel
D''-D\parallel= \parallel D'-D\parallel$ and such that
both $(D''-D)$ and $(D''-D')$ lie in the convex cone
spanned by the $A_i$. To see this, if
$D'-D=\sum_i\beta_iA_i$, simply take
$D''-D= \sum_i|\beta_i|A_i$.
\end{remark}

\noindent {\bf Proof of Theorem \ref{main}.} We start with two
lemmas. First, let $X$ be a smooth projective variety of dimension
$n$, and let $V\subset X$ be an irreducible subvariety of dimension
$d$. We recall the definition of the asymptotic order function
$\ord_V(\parallel \cdot \parallel)$ defined on ${\rm Big}(X)_{\RR}$
(we refer to \cite{elmnp} for the basic properties of this
function). If $D$ is a big $\QQ$-divisor, and if $m$ is divisible
enough, then $\ord_V(|mD|)$ denotes the order of vanishing at the
generic point of $V$ of a general element in $|mD|$. We have
$$\ord_V(\parallel
D\parallel):=\lim_{m\to\infty}\frac{\ord_V(|mD|)}{m}
=\inf_m\frac{\ord_V(|mD|)}{m}.$$ This extends as a continuous,
convex function to ${\rm Big}(X)_{\RR}$. Note that the notation in
the asymptotic order of vanishing of $D$ should not be confused with
the norm on $N^1(X)_{\RR}$.  For the rest of this section, we fix as
above a basis $A_1,\ldots,A_r$ for $N^1(X)_{\RR}$ consisting of
ample divisors, and let $\sigma$ be the cone generated by the $A_i$.
We consider the norm on $N^1(X)_{\RR}$ given by $\parallel\sum
\alpha_iA_i\parallel=\max_i|\alpha_i|$.

\begin{lemma}\label{order_growth} Given $\sigma$ as
above, there is $\beta>0$ such that for every big
$\RR$-divisor $D$ and for every $V\subseteq\B+(D)$, if
$A\in\sigma$ is nonzero and $D-A$ is big, we have
\begin{equation}\label{lower_bound}
\ord_V(\parallel D-A\parallel)\geq \beta\cdot\parallel
A\parallel.
\end{equation}
\end{lemma}

\begin{proof} We adapt to our more general setting the
argument for Lemma~1.4 in \cite{nakamaye2}. It is clear
by the continuity of the asymptotic order function that
it is enough to satisfy the condition in the statement
 for every nonzero
$A\in\sigma$ such that $D-A$ is a big $\QQ$-divisor.

Let $A'$ be a very ample divisor on $X$ such that
$T_X\otimes {\mathcal O}(A')$ is an ample vector
bundle. We can find a positive integer $b$ such that
for every nonzero $A\in\sigma$, the divisor
$\frac{b}{\parallel A\parallel}\cdot A-A'$ is ample. We
show now that $\beta=\frac{1}{b}$ satisfies our
requirement.  Suppose that there is a nonzero
$A\in\sigma$ such that $D-A$ is a big $\QQ$-divisor and
$\ord_V(\parallel D-A\parallel) < \beta\cdot\parallel
A\parallel$. This means that there is $m$ and
$s\in H^0(X,{\mathcal O}(m(D-A)))$ whose order of
vanishing at the generic point $\eta$ of $V$ is
$\ord_{\eta}(s)<m\beta\cdot\parallel A\parallel$. Note
that we may replace $m$ by any multiple, so we may
assume that $m$ is divisible enough.

We use the notation in \cite{eln}: if $B$ is a line
bundle on $X$, then ${\mathcal D}_B^{\ell}$ denotes the
\emph{bundle of differential operators} of order $\leq \ell$ on
$B$. This is defined as
$${\mathcal D}_B^{\ell} := \mathcal{H}om (P^{\ell}_B, B) =
{P^{\ell}_B}^{\vee} \otimes B, $$
where $P^{\ell}_B$ is the bundle of $\ell$-principal parts associated
to $B$ (having as fibers the spaces of $\ell$-jets of sections of $B$).
These bundles sit in short exact sequences of the form
$$0\longrightarrow  {\mathcal D}_B^{\ell - 1}\longrightarrow
{\mathcal D}_B^{\ell} \longrightarrow
{\rm Sym}^{\ell} (T_X)\longrightarrow 0.$$
It follows from Lemma~2.5 in \cite{eln} and our
hypothesis on $A'$ that there is $\ell_0$ such that for every $B$ and
every $\ell\geq\ell_0$, the sheaf ${\mathcal D}_B^{\ell}
\otimes{\mathcal O}_X(\ell A')$
is globally generated. Since every global section
of $B$ determines a vector bundle map ${\mathcal D}_B^{\ell} \rightarrow B$,
we have a natural induced map
$$H^0 (\mathcal{D}_B^{\ell}\otimes{\mathcal O}_X(\ell A')) \longrightarrow
H^0(\OO_X(B + \ell A^{\prime}))$$
which eventually produces
nontrivial sections in $H^0(\OO_X(B + \ell A^{\prime}))$
arising locally via the process of differentiation.

For our Lemma, by taking $B={\mathcal O}(m(D-A))$ it follows that if
$\ell\geq m\beta\cdot\parallel A\parallel$, then we may
apply a suitable differential operator to our section
$s$ to get $\widetilde{s}\in H^0(X, {\mathcal O}(m(D-A)+\ell A'))$
that does not vanish at $\eta$. If $mA-\ell A'$ is
ample, then this contradicts the assumption that $V$ is
contained in $\B+(D)$.

As $\frac{1}{\beta\cdot\parallel A\parallel}\cdot A-A'$
is ample by assumption, we may take $m$ large enough so
that
$$\frac{m}{m\beta\cdot\parallel A\parallel+1}\cdot
A-A'$$ is again ample. Moreover, we may assume that
$m\beta\cdot\parallel A\parallel\geq
\ell_0+1$. If we choose $\ell$ to be the smallest
integer $\geq m\beta\cdot\parallel A\parallel$, then
both our requirements on $\ell$ are satisfied. This
completes the proof of the lemma.
\end{proof}

Our next lemma deals with a subtracting procedure introduced in
\cite{nakamaye2} (note however that we add an extra
condition in order to fix a small gap in the proof in \emph{loc.
cit.}). The goal is to get a lower bound on the dimension of the
space of sections of $L$ minus an ample, starting from sections in
$L$ with small order of vanishing at given points.

We keep the assumption that $V$ is a subvariety of dimension $d\geq
1$ of the smooth, projective variety $X$, and let $L$ and $B$ be
$\QQ$-divisors on $X$ such that $B$ and $B-L$ are ample. We consider
integers $m$ and $k$ such that $m$ is divisible enough and $k\gg m$.
Suppose that we have positive rational numbers $\epsilon<\epsilon'$
such that if $x_1,\ldots,x_N\in V$ are general points, for every $m$
(divisible enough) $mL$ separates $m\epsilon'$ jets at
$x_1,\ldots,x_N$. Suppose that for every $m$ and $k$ we have a
closed subscheme $V_{m,k}$ of $X$ supported on $V$. If $I_V$ and
$I_{m,k}$ are the ideals defining $V$ and $V_{m,k}$, respectively,
we assume that around every smooth point of $V$ we have
$I_{m,k}\subseteq I_V^{km\epsilon}$.

Let us fix now $x_1,\ldots,x_N$ smooth points on $V$ as above such
that ${\rm depth}(\OO_{V_{m,k},x_i})=d$ for all $i$ (this assumption
is satisfied by general points). We fix also a positive rational
number $a$. Let $m_1$ be divisible enough such that $m_1aB$ is very
ample and the linear system
$$\Sigma=\{B_0\in |m_1aB|\mid x_i\in B_0\}$$
induces a basepoint-free linear system on
$X\smallsetminus\{x_1,\ldots,x_N\}$. If $s$ is a section of a line
bundle on a scheme $Z$, the order of vanishing of $s$ at a point
$x\in Z$ is the largest $p$ such that a local equation for $s$ at
$x$ lies in the $p$th power of the ideal defining $x$ in $Z$.

\begin{lemma}\label{subtraction}
With the above notation, suppose that for every $m$ and $k$ we have
a vector subspace $W_{m,k}\subseteq H^0(X,\OO(kmL))$ such that
\begin{enumerate}
\item[(i)] For every nonzero section $s$ in  $W_{m,k}$ we have ${\rm min}_i{\rm
ord}_{x_i}(s)< km\epsilon$.
\item[(ii)] If $B'$ is a very general divisor in $\Sigma$, then for every
$m$ and $k$, and every $s$ in $W_{m,k}$ such that $s\vert_{B'}$ is
nonzero, we have ${\rm min}_i{\rm ord}_{x_i}(s\vert_{B'})<
km\epsilon$.
\end{enumerate}
Then $W_{m,k}$ induces a vector subspace $W'_{m,k}\subseteq
H^0(V_{m,k},\OO(km(L-aB)))$ such that every nonzero section $s$ of
$W'_{m,k}$ satisfies
\begin{equation}\label{ineq100}
{\rm min}_i{\rm ord}_{x_i}(s)< km\epsilon
\end{equation}
 and such that
$$\dim\,W'_{m,k}\geq\dim\,
W_{m,k}-a\cdot\deg_B(V)\cdot\ell(\OO_{V_{m,k},\eta})\cdot(km)^d,$$
where $\eta$ denotes the generic point of $V$.
\end{lemma}

\begin{proof}
Note first that since around every $x_i$ we have $I_{m,k}\subseteq
I_V^{km\epsilon}$, a section of a line bundle on a subscheme $Z$ of
$X$ has order $< km\epsilon$ at $x_i$ if and only if this holds for
its restriction to $Z\cap V_{m,k}$.
 We assume that $m$ is divisible by $m_1$, and let
$B_1,\ldots,B_{km/m_1}$ be very general elements in the linear
system $\Sigma$, so they satisfy the condition (ii) above. Since
${\rm depth}(\OO_{V_{m,k},x_i})\geq 1$, we may also assume that no
$B_i$ contains an associated point of $V_{m,k}$.

Let $\overline{W}_{m,k}$ denote the image of $W_{m,k}$ in
$H^0(V_{m,k},\OO(kmL))$. Our assumption implies that the restriction
map gives an isomorphism $W_{m,k}\simeq\overline{W}_{m,k}$. Let
$W'_{m,k}$ be the kernel of the composition
$$\overline{W}_{m,k}\hookrightarrow H^0(V_{m,k}, \OO (kmL))
\to H^0(V_{m,k}, \OO_{V_{m,k}}(kmL)\vert_{\sum_iB_i}).$$ It is clear
that $W'_{m,k}\subseteq H^0(V_{m,k}, \OO(km(L-aB)))$ and that (i)
above implies that the nonzero sections in $W'_{m,k}$ satisfy
(\ref{ineq100}).

In order to get the lower bound for $\dim\,W'_{m,k}$, it is enough
to show that if
$$w_{p}:={\rm Im}\big{(} W_{m,k}\longrightarrow
H^0(V_{m,k}\cap B_{p},\OO_{V_{m,k}\cap B_{p}} (kmL))\big{)},$$ then
for every $p\leq km/m_1$ we have the following upper bound:
$$\dim(w_{p})\leq
m_1a\cdot\deg_B(V)\cdot\ell(\OO_{V_{m,k},\eta})\cdot(km)^{d-1}.$$

Since $B-L$ is ample, it follows that $mB$ also separates
$m\epsilon'$-jets at $x_1,\ldots,x_N$, hence we can find $m_2$ such
that if $D_1,\ldots,D_{d-1}$ are general elements in $|m_2B|$ with
order of vanishing $\geq m_2\epsilon+1$ at each $x_j$, then their
local equations
 form a regular sequence with respect to $V_{m,k}\cap
B_p$. We use here the fact that ${\rm depth}({\mathcal
O}_{V_{m,k}\cap B_p,x_j})=d-1$ for all $j$ and apply
Remark~\ref{extra_divisor} successively to avoid containing suitable
associated points.

 {}From now on we
assume that $m$ is divisible also by $m_2$. Property (ii) above
implies that no nonzero element in $w_i$ can lie in the image of
$$H^0\left(V_{m,k}\cap B_p, \OO\left(km\left(L-\frac{1}{m_2}\cdot
D_1\right)\right)\right) \longrightarrow H^0(V_{m,k}\cap
B_p,\OO(kmL)).$$ Therefore $\dim(w_p)$ is bounded above by the
dimension of the image of $W_{m,k}$ in $H^0(V_{m,k}\cap B_p\cap
km/m_2\cdot D_1, \OO(kmL))$. Moreover, since  $\frac{km}{m_2}D_1$
passes through $x_1,\ldots,x_N$ with multiplicity at least
$km\epsilon$, it follows that for every $s$ in $W_{m,k}$ such that
$s\vert_{B_p}$ is nonzero, we have
$${\rm min}_j{\rm ord}_{x_j}(s\vert_{B_p\cap D_1})< km\epsilon.$$
Therefore we can repeat the above procedure.
 After $d-1$ such steps, we deduce
$$\dim(w_p)\leq m_1a\cdot (km)^{d-1}\cdot
\deg_B(V_{m,k})=m_1a\cdot (km)^{d-1}\cdot\deg_B(V)
\cdot\ell(\OO_{V_{m,k},\eta}),$$ which completes the proof of the
lemma.
\end{proof}

\begin{remark}\label{moreover}
Suppose that in the above lemma we assume that all sections in
$W_{m,k}$ restrict to zero on subschemes $V'_{m,k}$ whose support is
properly contained in $V$. If the points $x_1,\ldots,x_N$ do not lie
in any support of a primary component of a $V'_{m,k}$ (which can be
achieved by taking the $x_i$ very general on $V$), then we get that
the sections in $W'_{m,k}$ also restrict to zero on $V'_{m,k}$.
Indeed, in the above proof it is enough to make sure that the
divisors $B_1,\ldots,B_{km/m_1}$ do not contain any of the supports
of the primary components of the schemes $V'_{m,k}$.
\end{remark}

We can give now the proof of Theorem~\ref{main}. We will use the
following notation: if $E$ is a divisor on $X$ such that
$|E|\neq \emptyset$, we will denote by $\frb_{|E|}$ the ideal
defining the base locus of this linear system.

\begin{proof}[Proof of Theorem~\ref{main}] The proof of
Theorem~\ref{main} follows the approach in
\cite{nakamaye2}, using in addition our result on
approximating volumes in terms of separation of jets.
Note that by Proposition~\ref{no_points}, we have
$\dim(V)=d\geq 1$. We may also assume that $D$ is big:
otherwise $V=X$, and the theorem follows from the
continuity of the usual volume function on
$N^1(X)_{\RR}$ (see \cite{positivity},
Corollary~2.2.45).  While the proof of the general case
is quite technical, if we assume that $V=\B+(D)$, then
the proof becomes more transparent. For the benefit of
the reader, we give first the proof of this particular
case, and we describe later the general argument.

We start therefore by assuming that $V=\B+(D)$ and that
$D$ contradicts the conclusion of the Theorem. It
follows from Remark~\ref{reduction} that
 there is $\delta>0$ and a sequence of ample divisors
$A_q$ going to $0$ and such that $D+A_q$ are $\QQ$-divisors with
$\vol_{X\vert V}(D+A_q)>\delta$ for every $q$. Moreover, we may
clearly assume that $V\not\subseteq\B+(D+A_q)$ for all $q$.

By Lemma~\ref{order_growth}, we can find $\beta\in\QQ_+^*$ such that
${\rm ord}_V(\parallel D-H\parallel)\geq\beta\cdot \parallel
H\parallel$ for every $H$ in the interior of our cone $\sigma$ such
that $D-H$ is big. In addition to this lower bound for the
asymptotic order function, we will also need an upper bound for the
asymptotic multiplicity. Recall from \cite{elmnp} that if $E$ is a
big $\QQ$-divisor such that $V$ is not properly contained in
$\BB(E)$, and if for $m$ divisible enough we denote by $e_m$ the
Samuel multiplicity of the local ring $\OO_{X,V}$ with respect to
the localization of $\frb_{|mE|}$, then this asymptotic multiplicity
is defined by
$$e_V(\parallel E\parallel):=\lim_{m\to\infty}\frac{e_m}{m^{n-d}}=
\inf_m\frac{e_m}{m^{n-d}}.$$ In our setup, since $V \not\subseteq
\B+(D+A_q)$ for any $q$, it follows that $e_V(\parallel
D\parallel)=0$. Moreover, $e_V(\parallel\cdot\parallel)^{1/(n-d)}$
is locally Lipschitz continuous, hence there is $M>0$ such that
\begin{equation}\label{bound_e}
e_V(\parallel D-H\parallel)< M\cdot \parallel H\parallel^{n-d}
\end{equation}
if $H$ lies in a suitable ball ${\mathcal U}$ around the origin. We
refer to \cite{elmnp}, the end of \S 2 and Remark~3.2 for the
properties of asymptotic multiplicity that we used (the fact that
$e_V(\parallel\cdot\parallel)^{1/(n-d)}$ is locally Lipschitz
continuous on its domain follows also from the fact that it is
homogeneous and convex).

We choose now a $\QQ$-divisor $B$ such that both $B$ and $B-D$ lie
in the interior of the cone $\sigma$. Let $a\in\QQ_+^*$ be small
enough, such that $\B+(D-aB)=\B+(D)$ and
\begin{equation}\label{ineq_a}
a<\frac{(n-d)!\cdot\delta\beta^{n-d}}{n!\cdot M\cdot\deg_B(V)}.
\end{equation}
We fix now $q$ such that $aB-A_q$ and $B-(D+A_q)$ lie in the
interior of $\sigma$ and from now on we put $A=A_q$.

Since $V\not\subseteq\B+(D+A)$ and $\vol_{X\vert V}(D+A)>\delta$,
Theorem~\ref{separation_of_jets} implies that there are $N\geq 1$
and $\epsilon\in\QQ_+^*$ such that $N\epsilon^d>\delta$ and if
$x_1,\ldots,x_N$ are general points on $V$, the canonical map
\begin{equation}\label{surj10}
H^0(X,\OO(p(D+A)))\to\bigoplus_{i=1}^N
H^0(X,\OO(p(D+A))\otimes\OO_{X,x_i}/\frmm_{x_i}^{p\epsilon})
\end{equation}
 is surjective for $p$ divisible enough.
 Moreover, by taking $N$ large enough we can make sure that
 $\epsilon$ is as small as we want.

We fix now $H$ small enough in the interior of $\sigma\cap {\mathcal
U}$ such that $aB-A-H$ is an ample $\QQ$-divisor. Let us choose
$\epsilon$ and $N$ as above such that $\epsilon/\beta<\parallel
H\parallel$. After subtracting from $H$ a small multiple of a
$\QQ$-divisor, we can make $\parallel H\parallel$ arbitrarily close
to $\epsilon/\beta$, and therefore by (\ref{ineq_a}) we may assume
that
\begin{equation}\label{ineq_a2}
a<\frac{(\epsilon/\beta)^{n-d}}{\parallel H\parallel^{n-d}}\cdot
\frac{(n-d)!\cdot\delta\beta^{n-d}}{n!\cdot M\cdot\deg_B(V)}.
\end{equation}
We see that $D-H$ is a $\QQ$-divisor and
$\B+(D)\subseteq\BB(D-H)\subseteq\B+(D-aB)$, hence $\BB(D-H)=V$.
Note that we may assume in addition that (\ref{surj10}) is
surjective also for some $\epsilon'>\epsilon$, if $p$ is divisible
enough.

We  will consider integers $m$ and $k$, with $m$ large and divisible
enough, and with $k\gg m$. For every such $k$ and $m$, let
$\pi_m\colon X_m\to X$ be a log resolution of $\frb_{|m(D-H)|}$, and
write $(\pi_m)^{-1}(|m(D-H)|)=E_m+|M_m|$, with $E_m$ the fixed part
and $M_m$ the moving part. The ideal $I_{m,k}=(\pi_m)_*\OO_{X_m}(-kE_m)$ is
equal to the integral closure of $\frb_{|m(D-H)|}^k$ and it defines
a subscheme that we denote by $V_{m,k}$. Note that since
$\BB(D-H)=V$ and $m$ is divisible enough, $V_{m,k}$ is supported on
$V$. We denote by $I_V$ the ideal defining $V$. By our choice of
$\beta$, and since $\epsilon/\beta<\parallel H\parallel$, we see
that in the neighborhood of any smooth point of $V$ we have
\begin{equation}
I_{m,k}\subseteq I_V^{km\epsilon}
\end{equation}
(we use the fact that for an ideal defining a smooth subvariety, the
symbolic powers coincide with the usual powers and they are
integrally closed).

Let $m_1$ be divisible enough, so $m_1aB$ is very ample, the linear
system
$$\Sigma=\{B_0\in |m_1aB|\vert x_i\in B_0\,{\rm for}\,{\rm
all}\,i\}$$
 has no base points in
$X\smallsetminus\{x_1,\ldots,x_N\}$, and a general element of
$\Sigma$ is smooth. We choose
 $N$ general points $x_1,\ldots,x_N$ such that the above properties
 of $\Sigma$ are satisfied, and in addition
 (\ref{surj10}) is
surjective and ${\rm depth}(\OO_{V_{m,k},x_i})=d$ for every $i$.
Moreover, we can choose these points such that if $p$ is divisible
enough, then $p(D+A)-m_1aB$ separates $p\epsilon$-jets at
$x_1,\ldots,x_N$ (see Remark~\ref{separation3}).

Our plan now is to apply Lemma~\ref{subtraction}. To this end, fix
any very general divisor $B_0\in |m_1aB|$ passing through
$x_1,\ldots,x_N$. We have the following commutative diagram

\[
\begin{CD}
H^0(X,\OO(km(D+A)-m_1aB)))@>{\alpha_{m,k}}>>H^0(X,\OO(km(D+A)))\\
@VV{\phi_{m,k}}V @VV{\psi_{m,k}}V\\
\oplus_{i=1}^NH^0(\OO(km(D+A)-m_1aB)\otimes\OO_{X,x_i}
/\frmm_{x_i}^{km\epsilon})@>{\beta_{m,k}}>>\oplus_{i=1}^NH^0(\OO(km(D+A))\otimes\OO_{X,x_i}
/\frmm_{x_i}^{km\epsilon})
\end{CD}
\]
where the horizontal maps are induced by local equations of $B_0$.
Note that $\alpha_{m,k}$ is injective, and by construction,
$\phi_{m,k}$ and $\psi_{m,k}$ are surjective.
 Therefore we can choose
$W_{m,k}\subseteq H^0(X,\OO(km(D+A)))$ such that $W_{m,k}$ is mapped
isomorphically by $\psi_{m,k}$ onto
 $\oplus_{i=0}^NH^0(X,\OO(km(D+A))\otimes\OO_{X,x_i}/\frmm_{x_i}^{km\epsilon})$
and such that ${\rm Im}(\beta_{m,k})$ is in the image of
$W_{m,k}\cap {\rm Im}(\alpha_{m,k})$. This implies that for every
nonzero section $s$ in $W_{m,k}$, we have
\begin{equation}
\min_i{\rm ord}_{x_i}(s)< km\epsilon.
\end{equation}
 Moreover, if $s\vert_{B_0}$ is nonzero, then
\begin{equation}\label{eq200}
\min_i{\rm ord}_{x_i}(s\vert_{B_0})< km\epsilon.
\end{equation}

We assert moreover that the analogue of (\ref{eq200}) holds for any
very general $B'_0$ in $|m_1aB|$ passing through the $x_i$. Indeed,
consider for such $B'_0$ the commutative diagram
\[
\begin{CD}
W_{m,k}@>>>H^0(B'_0,\OO(km(D+A)))\\
@VV{\psi_{m,k}}V @VV{\rho_{m,k}}(B'_0)V\\
J_{m,k}@>>>\overline{J}_{m,k}(B'_0)
\end{CD}
\]
where $J_{m,k}=\bigoplus_{i=1}^N H^0(X,\OO(km(D+A))\otimes
\OO_{X,x_i}/\frmm_{x_i}^{km\epsilon})$ and
$$\overline{J}_{m,k}(B'_0)= \bigoplus_{i=1}^N H^0(B'_0,
\OO(km(D+A))\otimes\OO_{B'_0,x_i}/\frmm_{x_i}^{km\epsilon}).$$ Since
$\psi_{m,k}$ is surjective, we see that the restriction of
$\rho_{m,k}(B'_0)$ to the image of $W_{m,k}$ in
$H^0(B'_0,\OO(km(D+A)))$ is always surjective. Our assertion is that
this image maps isomorphically onto $\overline{J}_{m,k}(B'_0)$. By
construction, this holds when $B'_0=B_0$, and since the spaces
involved have constant dimension for general $B'_0$, it follows that
(\ref{eq200}) holds with $B_0$ replaced by $B'_0$. Therefore the two
hypotheses (i) and (ii) in Lemma~\ref{subtraction} are satisfied for
$L=D+A$.

By construction we have
\begin{equation}
\dim\,W_{m,k}=\sum_{i=1}^N\dim\,\OO_{X,x_i}/\frmm_{x_i}^{km\epsilon}
=N{{km\epsilon+n-1}\choose{n}}.
\end{equation}
Since $N \epsilon^d > \delta$, we deduce
\begin{equation}\label{bound_below2}
\dim\,W_{m,k}
>\frac{\delta\epsilon^{n-d}}{n!}(km)^n+O((km)^{n-1}).
\end{equation}

On the other hand, Lemma~\ref{subtraction} gives a vector subspace
$$W'_{m,k}\subseteq H^0(V_{m,k},\OO(km(D+A-aB)))$$
such that
\begin{equation}\label{bound_below}
\dim\,W'_{m,k}\geq\dim\,W_{m,k}-a\cdot\deg_B(V)\cdot
\ell(\OO_{V_{m,k},\eta})\cdot (km)^d
\end{equation}
and such that for every nonzero section $s$ in $W'_{m,k}$, we have
\begin{equation}\label{bound_ord}
\min_{i=1}^N{\rm ord}_{x_i}(s)< km\epsilon.
\end{equation}
Since $aB-A-H$ is ample, we get corresponding spaces of sections
$$W''_{m,k}\subseteq H^0(V_{m,k},\OO(km(D-H)))$$
satisfying the same lower bound on the dimension and such that for every
nonzero section in $W''_{m,k}$ we have (\ref{bound_ord}).

 We give now an upper bound for
$\ell(\OO_{V_{m,k},\eta})$ when $m$ is divisible enough, but fixed,
and $k$ goes to infinity. We clearly have
\begin{equation}
\ell(\OO_{V_{m,k},\eta})\leq \ell(\OO_{X,\eta}/\frb_{|m(D-H)|}^k)<
\frac{\tilde{e}_m}{(n-d)!}\cdot k^{n-d}
\end{equation}
if $\tilde{e}_m$ is larger than the multiplicity $e_m$ of
$\OO_{X,\eta}$ with respect to the localization of
$\frb_{|m(D-H)|}$. By (\ref{bound_e}), since $m$ is large enough we
may choose such $\tilde{e}_m$ with $\tilde{e}_m<M\cdot
m^{n-d}\parallel H\parallel^{n-d}$ and conclude that
\begin{equation}\label{bound_for_l}
\ell(\OO_{V_{m,k},\eta})< \frac{M\cdot\parallel
H\parallel^{n-d}}{(n-d)!}(km)^{n-d}\,\,{\rm for}\,k\gg 0.
\end{equation}
 Combining this with
(\ref{bound_below}), (\ref{bound_below2}) and (\ref{ineq_a2}) we
deduce that if $m$ is large and divisible enough, then
$\dim\,W''_{m,k}$ grows like a polynomial of degree $n$ in $k$, when
$k$ goes to infinity.

We use this to construct global sections of $km(D-H)$. Consider the
exact sequence
$$H^0(X,\OO(km(D-H)))\longrightarrow
H^0(V_{m,k},\OO(km(D-H)))\longrightarrow
H^1(X,\OO(km(D-H))\otimes I_{m,k}).$$
Since $M_m$
is nef, we have $h^1(X_m,\OO(kM_m))\leq O(k^{n-1})$ for $k\gg 0$.
Using the Leray spectral sequence, this gives
$$h^1(X, \OO_X(km(D-H))\otimes
I_{m,k})\leq O(k^{n-1})$$ for $k\gg 0$. Therefore
most of the sections in $W''_{m,k}$ can be lifted to
$H^0(X,\OO(km(D-H)))$.

On the other hand, recall that for every nonzero section $s$ in
$W''_{m,k}$ we can find a point $x_i$ such that ${\rm ord}_{x_i}(s)<
km\epsilon$. Since ${\rm ord}_V(|km(D-H)|)\geq km\beta\parallel
H\parallel$ and $\epsilon/\beta<\parallel H\parallel$, we get a
contradiction. This completes the proof in the case $\B+(D)=V$.

\bigskip

We treat now the general case. The above proof fails since the
subschemes $V_{m,k}$ as defined above are not supported on $V$
anymore. We need to do some extra work to ensure that the sections
we construct on $V_{m,k}$ can be extended to subschemes supported on
the whole $\B+(D)$. We will use the following notation: if $F$ is an
integral divisor such that $|F|\neq\varnothing$, then we denote the
integral closure of $\frb_{|F|}^k$ by $\frb_{|F|}^{(k)}$. We
consider also $\frb_{|F|}^{(k)_2}$ that is defined locally as the
ideal of sections $\phi$ in $\OO_X$ such that for every divisor $T$
over $X$, with center on $X$ different from $V$, we have ${\rm
ord}_T(\phi) \geq k\cdot {\rm ord}_T(\frb_{|F|})$. Note that for
every $F_1$ and $F_2$, it follows from definition that
\begin{equation}\label{inclusion_of_ideals}
\frb_{|mF_1|}^{(k)_2}\cdot\frb_{|mF_2|}^{(k)_2}\subseteq\frb
_{|mF_1+mF_2|}
^{(k)_2}.
\end{equation}
It is clear that if $V$ is not contained in $\BB(F)$, then
$\frb_{|mF|}^{(k)}=\frb_{|mF|}^{(k)_2}$ for $m$ divisible enough. On
the other hand, if $V$ is an irreducible component of $\BB(F)$ and
$m$ is divisible enough, then $\frb_{|mF|}^{(k)}$ has a uniquely
determined primary component supported on $V$. If this is defined by
$\frb_{|mF|}^{(k)_1}$, then we have $\frb_{|mF|}^{(k)}=
\frb_{|mF|}^{(k)_1}\cap\frb_{|mF|}^{(k)_2}$.

We assume for the moment the following technical lemma.

\begin{lemma}\label{ideal_I} Fix an ample divisor $\overline{B}$
such that $D-\overline{B}$ is a $\QQ$-divisor with
$\BB(D-\overline{B})=\B+(D)$. There is a sheaf of ideals ${\mathcal
I}$ on $X$ whose support does not contain $V$ such that
$${\mathcal I}^{kp}\cdot
\frb_{|p\gamma(D+A)|}^{(k)}
\subseteq\frb_{|p(\gamma+1)(D+A-\frac{1}{\gamma}B_2)|}^{(k)_2}$$ for
every $A$, $B_2$, $\gamma$, $p$ and $k$ as follows: $A$ is ample
such that $D+A$ is a $\QQ$-divisor with $V\not\subseteq \BB(D+A)$,
$B_2$ is an ample $\QQ$-divisor such that $\overline{B}-2B_2$ is
ample, $\gamma>1$ a rational number,  $p$ is divisible enough
(depending on $A$, $B_2$, and $\gamma$), and $k$ is an arbitrary
positive integer.
\end{lemma}

In order to prove the general case of the theorem, we start by
fixing $\delta$, $\beta$ and $M$ as before. We fix also
$\overline{B}$ and ${\mathcal I}$ as in the above lemma. Let $B_1$
be an integral divisor  such that both $B_1$ and $B_1-D$ lie in the
interior of $\sigma$ and ${\mathcal I}\otimes\OO(B_1)$ is globally
generated. Let $B_2$ be a divisor in the interior of $\sigma$ such
that $\overline{B}-2B_2$ is ample and we put $B=B_1+B_2$.

We choose now $a\in\QQ_+^*$ small enough such that
$\B+(D-aB)=\B+(D)$, $a<1$ and
\begin{equation}\label{ineq_a3}
a(1+a)^d<\frac{(n-d)!\cdot\delta\beta^{n-d}}{n!\cdot
M\cdot\deg_B(V)}.
\end{equation}

As before, we can choose $A=A_q$ for $q\gg 0$ such that $aB_2-A$
lies in the interior of $\sigma$. Moreover, we can choose
$\epsilon\in\QQ_+^*$ and $N\geq 1$ such that $N\epsilon^d>\delta$
and (\ref{surj10}) is surjective for very general points
$x_1,\ldots,x_N$ on $V$ if $p$ is divisible enough (we may assume
that the same map is surjective also for some $\epsilon'>\epsilon$).
Arguing as before, we can find a divisor $H$ in the interior of
$\sigma\cap {\mathcal U}$ such that $aB_2-A-H$ is an ample
$\QQ$-divisor and we have $\epsilon/\beta<\parallel H\parallel$ and
\begin{equation}\label{ineq_a4}
a(a+1)^d<\frac{(\epsilon/\beta)^{n-d}}{\parallel
H\parallel^{n-d}}\cdot \frac{(n-d)!\cdot\delta\beta^{n-d}}{n!\cdot
M\cdot\deg_B(V)}.
\end{equation}

Note that $D-H$ is a $\QQ$-divisor with $\BB(D-H)=\B+(D)$. Let
$V_{m,k}$ be the subscheme defined by the primary ideal
$\frb_{|m(D-H)|}^{(k)_1}$, so $V_{m,k}$ has support $V$. We choose
again $m_1$ divisible enough, and very general points
$x_1,\ldots,x_N$ on $V$ such that $m_1aB$ is  very ample, the linear
subsystem of $|m_1aB|$ passing through $x_1,\ldots,x_N$  has no base
points in $X\smallsetminus\{x_1,\ldots,x_N\}$ and a general element
is smooth, (\ref{surj10}) is surjective and ${\rm
depth}(\OO_{V_{m,k},x_i})=d$ for every $i$. Moreover, if $p$ is
divisible enough, then $p(D+A)-m_1aB$ separates $p\epsilon$-jets at
$x_1,\ldots,x_N$. We may assume also that the points $x_i$ do not
lie on the support of any irreducible component of the schemes
defined by the ideals $\frb_{|m(D+A)|}^{(k)}$.

 We claim that we can find vector subspaces
$$\widetilde{W}_{m,k}\subseteq H^0\left(X,\OO(km(a+1)(D+A))\otimes \frb^{(k)}
_{|m(D+A)|}\right)$$ such that the induced map
\begin{equation}\label{eq101}
\widetilde{W}
_{m,k}\to \bigoplus_{i=1}^N
H^0(X,\OO(km(a+1)(D+A))\otimes\OO_{X,x_i}/\frmm_{x_i}^{km\epsilon+1})
\end{equation}
is an isomorphism and for a very general $B_0\in |m_1aB|$ passing
through $x_1,\ldots,x_N$ and for every $s\in \widetilde{W}_{m,k}$
with $s\vert_{B_0}\neq 0$, there is a point $x_i$ such that ${\rm
ord}_{x_i}(s\vert_{B_0})\leq km\epsilon$. Indeed, suppose first that
$B_0$ is very general as above, but fixed. Arguing as before, we see
that we can find $W_{m,1}\subseteq H^0(X,\OO(m(D+A)))$ that maps
isomorphically onto
$$\bigoplus_{i=1}^NH^0(X,\OO(m(D+A))\otimes\OO_{X,x_i}/\frmm_{x_i}^{m\epsilon})$$
and a subspace $W_{m,1}^{\circ}\subseteq W_{m,1}$ that maps
isomorphically onto
$$\bigoplus_{i=1}^NH^0(X,\OO(m(D+A)-m_1aB_0)\otimes\OO_{X,x_i}/\frmm_{x_i}^{m\epsilon}).$$
Note also that we have
 $H^0(X,\OO(m(D+A)))=H^0(X,\OO(m(D+A))\otimes\frb_{|m(D+A)|})$.

We can choose now subspaces $W_{m,k}^{\circ}\subseteq
W_{m,k}\subseteq H^0(X,\OO(km(D+A)))$ such that $W_{m,k}$ maps
isomorphically onto
$$\bigoplus_{i=1}^NH^0(X,\OO(km(D+A))\otimes\OO_{X,x_i}/\frmm_{x_i}^{km\epsilon}),$$
$W^{\circ}_{m,k}$ maps isomorphically onto
$$\bigoplus_{i=1}^NH^0(X,\OO(km(D+A)-m_1aB_0)\otimes\OO_{X,x_i}/\frmm_{x_i}^{km\epsilon})$$
and $W_{m,k}$ is contained in the image of
$${\rm Sym}^kW_{m,1}\to H^0(X,\OO(km(D+A))),$$
while $W_{m,k}^{\circ}$ is contained in the image of
$$W_{m,1}^{\circ}\otimes{\rm Sym}^{k-1}W_{m,1}\to H^0(X,\OO(km(D+A))).$$
Note that $W_{m,k}^{\circ}$ is the analogue of the image of the map
$\alpha_{m,k}$ in the top diagram on p.26.

 By construction, it follows that $W_{m,k}\subseteq H^0(X,
\OO(km(D+A))\otimes \frb^{(k)}_{|m(D+A)|})$. On the other hand,
since $V$ is not contained in $\BB(D+A)$ and $m$ is divisible enough
(recall also that the points $x_i$ are very general on $V$), there
are sections $t_{m,k}\in H^0(X,\OO(kma(D+A)))$ that do not vanish at
any of the points $x_i$. Multiplying by the section $t_{m,k}$
induces an embedding of $W_{m,k}$ in $H^0(X,\OO(km(a+1)(D+A))\otimes
\frb_{|m(D+A)|}^k)$. If we denote by $\widetilde{W}_{m,k}$ its
image, then it is clear that it satisfies the claimed properties. We
deduce as before that if we replace $B_0$ by a very general element
$B'_0\in |m_1aB|$ passing through $x_1,\ldots,x_N$, then it is still
true that every nonzero restriction to $B'_0$ of an element in
$\widetilde{W}_{m,k}$ has order $<km\epsilon$ at some $x_i$.

In particular, the $\widetilde{W}_{m,k}$ also satisfy the lower
bound (\ref{bound_below2}). We apply Lemma~\ref{subtraction} and
Remark~\ref{moreover} for $L=(D+A)$, with $m(a+1)$ instead of $m$,
to get as before spaces of sections
$$W'_{m,k}\subseteq H^0\left(V_{m,k},\OO(km(a+1)(D+A-aB))\otimes
\frb_{|m(D+A)|}^{(k)}\right)$$ such that
$$\dim\,W'_{m,k}\geq\dim\,\widetilde{W}_{m,k}-
a(a+1)^d\cdot\deg_B(V)\cdot\ell(\OO_{V_{m,k},\eta})\cdot (km)^d.$$
On the other hand, the lower bound (\ref{bound_for_l}) still holds,
and combining this with (\ref{bound_below2}) and (\ref{ineq_a4}) we
deduce that if $m$ is divisible enough, then $\dim\,W'_{m,k}$ grows
like a polynomial of degree $n$ in $k$ when $k$ goes to infinity.

Now since ${\mathcal I}\otimes\OO(B_1)$ is globally generated,
the support of ${\mathcal I}$ does not contain $V$, and $V_{m,k}$
is defined by a primary ideal,
we get an embedding
\begin{equation}\label{emb1}
W'_{m,k}\hookrightarrow
H^0\left(V_{m,k},\OO(km(a+1)(D+A-aB_2))\otimes {\mathcal
I}^{kma(a+1)} \frb_{|m(D+A)|}^{(k)}\right).
\end{equation}
Since ${\mathcal I}^{kma(a+1)}\subseteq {\mathcal I}^{kma}$,
applying Lemma~\ref{ideal_I} for $p=ma$ and $\gamma=1/a$, we deduce
\begin{equation}\label{emb2}
W'_{m,k}\hookrightarrow H^0\left(V_{m,k}, \OO(km(a+1)(D+A-aB_2))
\otimes\frb_{|m(a+1)(D+A-aB_2)|}^{(k)_2}\right).
\end{equation}
Furthermore, since $aB_2-A-H$ is ample, we get
an embedding
\begin{equation}\label{emb3}
W'_{m,k}\hookrightarrow H^0\left(V_{m,k}, \OO(km(a+1)(D-H))\otimes
\frb_{|m(a+1)(D-H)|}^{(k)_2}\right).
\end{equation}
Note that we may assume that (\ref{emb1}) and (\ref{emb3}) are
induced by multiplication with sections that do not vanish at any of
the $x_i$, so they do not increase the order of vanishing at these
points.

Because all sections in $W'_{m,k}$ vanish on the subscheme defined
by $\frb_{|m(D-H)|}^{(k)_2}$, it follows that they can be extended
by zero to give a space $W''_{m,k}$ of sections of
$\OO(km(a+1)(D-H))$ on the subscheme defined by
$\frb_{|m(D-H)|}^{(k)}$. As before, we show using a resolution
$\pi_m\colon X_m\to X$ that for $k\gg 0$, most of the sections in
$W''_{m,k}$ extend to $X$. On the other hand, for every nonzero
section $s$ in $W''_{m,k}$ we can find $i$ such that ${\rm
ord}_{x_i}(s)\leq km \epsilon-1$. Since by our choice of $\beta$ we
have
$${\rm ord}_V(|km(a+1)(D-H)|)\geq\beta km(a+1)\parallel H\parallel
>km\epsilon(a+1),$$
this gives a contradiction and completes the proof of the
Theorem in the general case.

\end{proof}

\begin{proof}[Proof of Lemma~\ref{ideal_I}] As
$\frb_{|p\gamma(D+A)|}^{(k)}=\frb_{|p\gamma(D+A)|}^{(k)_2}$,
it follows from (\ref{inclusion_of_ideals}) that it is
enough to have
$${\mathcal
I}^{kp}\subseteq\frb_{|p(D+A-\frac{\gamma+1}{\gamma}B_2)|}
^{(k)_2}.$$ Moreover, if we have this inclusion for $k=1$, then we
get it for all $k$.  Note that
$(D+A-\frac{\gamma+1}{\gamma}B_2)-(D-\overline{B})$ is ample, so
that $\frb_{|p(D-\overline{B})|}^{(1)_2}\subseteq
\frb_{|p(D+A-\frac{\gamma+1}{\gamma}B_2)|}^{(1)_2}.$ Therefore it is
enough to choose $q$ such that $q(D-\overline{B})$ is integral and
the reduced base locus of $|q(D-\overline{B})|$ is $\B+(D)$, and to
take ${\mathcal I}=\frb_{|q(D-\overline{B})|}^{(1)_2}$.
\end{proof}

\section{Moving Seshadri constants}

Moving Seshadri constants have been introduced in
\cite{nakamaye2} for the description of the augmented
base locus. In the case of nef line bundles, they
coincide with the usual Seshadri constants. In this
section we prove the basic properties of these
invariants, and we use the results in the previous
sections to deduce a stronger version of the main
result in \cite{nakamaye2}.  If $D$ is a nef
$\QQ$-divisor on a smooth, projective variety $X$, we
denote by $\epsilon(D;x)$ the Seshadri constant of $D$
at $x$. For the definition and basic results on
Seshadri constants we refer to
\cite{positivity} \S 5.1.

As in the case of asymptotic intersection numbers,
there are two equivalent definitions for moving
Seshadri constants. We start this time with the
definition in terms of arbitrary decompositions for the
pull-back of $D$, a definition which applies to
arbitrary $\RR$-divisors (note the similarity with the
formula in Proposition~\ref{general_resolution}).
Suppose that
$x\in X$ and that $D$ is a divisor such that
$x\not\in\B+(
D)$. We consider projective morphisms $f :
X'\longrightarrow X$, with $X'$ smooth, which are
isomorphisms over a neighborhood of $x$, and
decompositions $f^*(D)=A+E$, with $A$ an ample
$\QQ$-divisor and $E$ effective such that $f^{-1}(x)$
is not in the support of $E$. Note that for every such
$f$, we have $f^{-1}(x)\not\in\B+(f^*(D))$, so there
exist indeed decompositions as described above.
\begin{definition}\label{definition_Seshadri} Let $D$
be an $\RR$-divisor. If $x\not\in\B+(D)$, then the
\emph{moving Seshadri constant} of $D$ at $x$ is
\begin{equation}
\epsilon(\parallel D\parallel;
x):=\sup_{f^*(D)=A+E}\epsilon(A,x),
\end{equation} where the supremum is over all morphisms
$f$ and decompositions
$f^*(D)=A+E$ as above. If $x \in \B+(D)$, then we put
$\epsilon(\parallel D\parallel;x)=0$.
\end{definition}

It is easy to see that the above invariant is finite
(see, for example, Proposition~\ref{basic1} ${\rm i)}$
below). It is also clear from the definition that $x$
is in $\B+(D)$ if and only if
$\epsilon(\parallel D\parallel;x)=0$. The value $0$
over $\B+$ is justified by the following theorem, which
is our main result on moving Seshadri constants. As we
will see, it can be considered a stronger version of
Theorem~0.8 in \cite{nakamaye2}.
\begin{theorem}\label{continuity_Seshadri} For every
point $x$ in $X$, the map
$D\longrightarrow\epsilon(\parallel D\parallel;x)$ is
continuous on the entire N\'{e}ron-Severi space
$N^1(X)_{\RR}$.
\end{theorem}

The proof will be given at the end of this section. We
start by giving some basic properties and
interpretations of the moving Seshadri constants. As
moving Seshadri constants of non-big divisors are
trivial, we henceforth assume that all divisors are
big.
\begin{proposition}\label{basic1} Suppose that $D$ is a
big $\RR$-divisor on $X$.
\item{\rm i)} We have $\epsilon(\parallel D\parallel;x)
\leq \vol_X(D)^{1/\dim(X)}$.
\item{\rm ii)} If $D\equiv E$, then $\epsilon(\parallel
D\parallel;x) =\epsilon(\parallel E\parallel;x)$.
\item{\rm iii)} $\epsilon(\parallel \lambda
D\parallel;x)=
\lambda\cdot\epsilon(\parallel D\parallel;x)$ for every
positive $\lambda$.
\item{\rm iv)} If $D$ is an ample $\QQ$-divisor, then
$\epsilon(\parallel D\parallel;x)=\epsilon(D;x)$.
\item{\rm v)} If $D'$ is another $\RR$-divisor such
that $x\notin \B+(D) \cup \B+(D')$, then
\begin{equation}
\epsilon(\parallel D+D'\parallel;x)\geq
\epsilon(\parallel D\parallel;x)+\epsilon(\parallel
D'\parallel;x).
\end{equation}
\end{proposition}

\begin{proof} All proofs follow from definition and
from the properties of the usual Seshadri constants.
\end{proof}

We explain now the connection with the definition of
moving Seshadri constants from \cite{nakamaye2}. This
is analogous to the definition of asymptotic
intersection numbers. Suppose that $D$ is a
$\QQ$-divisor and that $x \notin \BB(D)$. Let $m$ be
sufficiently divisible, so $mD$ is an integral divisor
and $x$ is not in the base locus of $|mD|$. We take a
resolution of this base locus as in
Definition~\ref{Asymptotic intersection number} (with
$V$ replaced by $x$).

We define following \cite{nakamaye2}:
\begin{equation}
\epsilon'(\parallel D\parallel;x):
=\lim_{m\to\infty}\frac{\epsilon(M_m;\pi_m^{-1}(x))}
{m}=\sup_m\frac{\epsilon(M_m;\pi_m^{-1}(x))}{m},
\end{equation} where the limit and the supremum are
over divisible enough $m$. Note that
$\epsilon(M_m;x)$ does not depend on the particular
morphism $\pi_m$. Moreover, given positive integers $p$
and $q$, sufficiently divisible , we may choose
$\pi : X'\longrightarrow X$ that satisfies our
requirements for
$|pD|$, $|qD|$ and $|(p+q)D|$. Since we have
$M_{p+q}=M_p+M_q+E$, for some effective divisor $E$ with
$\pi^{-1}(x)\not\in {\rm Supp}(E)$, we deduce  that
$$\epsilon(M_{p+q};\pi^{-1}(x))\geq\epsilon(M_p;\pi^{-1}(x))
+\epsilon(M_q,\pi^{-1}(x)).$$ This implies that the
limit in the definition of $\epsilon'(\parallel
D\parallel;x)$ exists, and it is equal to the
corresponding supremum.  We now show that the two
invariants we have defined are the same.

\begin{proposition}\label{numerical_description} If $D
$
is a big $\QQ$-divisor and if $x\not\in \BB(D)$, then
$\epsilon'(\parallel D\parallel;x)=\epsilon(\parallel
D\parallel;x)$.
\end{proposition}

\begin{proof} By replacing $D$ with a suitable
multiple, we may assume that $D$ is integral,
$\BB(D)={\rm Bs}(|D|)_{\rm red}$ and that $|D|$ defines
a rational map whose image has dimension $n$. For
$m\in\NN^*$, take $\pi_m : X_m\longrightarrow X$ as in
the definition of $\epsilon'(\parallel D\parallel;x)$
and write
$\pi_m^*(mD)=M_m+E_m$. Recall our assumption that
$\pi_m^{-1}(x)\not\in {\rm Supp}(E_m)$.

Suppose first that $x\in \B+(D)$. In this case, it
follows
 easily that $\pi_m^{-1}(x)\in \B+(M_m)$. Since $M_m$
is big and nef and $\pi_m^{-1}(x)\in \B+(M_m)$,
Corollary~\ref{cor_nakamaye} implies that there is a
subvariety
$V\subseteq X_m$ of dimension $d\geq 1$, such that
$\pi_m^{-1}(x)\in V$  and $(M_m^d\cdot V)=0$. Therefore
$\epsilon(M_m;\pi_m^{-1}(x)) =0$, and since this is
true for every $m$ we get $\epsilon'(\parallel
D\parallel;x)=0$.

Suppose now that $x\notin\B+(D)$.  We show first that
$\epsilon(\parallel
D\parallel;x)\geq\epsilon'(\parallel D\parallel;x)$ by
proving that for every $m$ divisible enough,
$\epsilon(\parallel D\parallel;x)\geq
\epsilon(M_m;\pi_m^{-1}(x))/m$. If $\pi_m^{-1}(x)\in
\B+(M_m)$, then the above argument using
Corollary~\ref{cor_nakamaye} shows that
$\epsilon(M_m;\pi_m^{-1}(x))=0$, and we are done.

Therefore we may assume that
$\pi_m^{-1}(x) \notin \B+(M_m)$, so we can write
$M_m=A+E$, where $A$ is ample, $E$ is effective, and
$\pi_m^{-1}(x)\not\in{\rm Supp}(E)$. If $p\in\NN^*$,
then we have
$M_m=(1/p)E+A_p$, where
$A_p=\frac{1}{p}A+\frac{p-1}{p}M_m$ is ample. It
follows from definition that $\epsilon(\parallel
D\parallel;x)
\geq(1/m)\epsilon(A_p;\pi_m^{-1}(x))$ for every $p$. By
letting $p$ go to infinity, we deduce
$\epsilon(\parallel D\parallel;x)
\geq\epsilon(M_m;\pi_m^{-1}(x))/m$.

We prove now that $\epsilon(\parallel D\parallel;x)\leq
\epsilon'(\parallel D\parallel;x)$. Let
$f:X'\longrightarrow X$ and $f^*(D)=A+E$ be as in t
he
definition of $\epsilon(\parallel D\parallel;x)$. Fix
$m$ such that $mA$ is integral and very ample. By
taking a log resolution of the base locus of $f^*(mD)$
which is an isomorphism over a neighborhood of
$f^{-1}(x)$, we may assume that we can write $f^*(mD)=
M_m + E_m$ as in the definition of
$\epsilon'(\parallel D\parallel;x)$. Since $mA$ is
basepoint-free, we have
$M_m=mA+E_m'$, where $E'_m$ is effective and
$E_m'\leq mE$, so $f^{-1}(x) \notin {\rm Supp}(E_m')$.
Therefore $\epsilon(M_m;f^{-1}(x))/m\geq
\epsilon(A;f^{-1}(x))$, hence $\epsilon'(\parallel
D\parallel;x)\geq\epsilon(\parallel D\parallel;x)$, and
this completes the proof of the Proposition.
\end{proof}

\begin{remark}\label{nef_Seshadri} If $D$ is a
$\QQ$-divisor that is nef and big, then
$\epsilon(\parallel D\parallel;x)=\epsilon(D;x)$. This
follows of course from the corresponding property for
ample divisors, together with the continuity property
of both invariants (see
Theorem~\ref{continuity_Seshadri} above). However, we
can also give a direct argument as follows.  If
$x\in\B+(D)$, then $\epsilon(\parallel D\parallel;x)=0$
by definition, while $\epsilon(D;x)=0$ by
Corollary~\ref{cor_nakamaye}. Suppose now that $x
\notin \B+(D)$. If $f : X'\longrightarrow X$ and
$f^*(D)=A+E$ are as in the definition of
$\epsilon(\parallel D\parallel;x)$, using the fact that
 $x$ is not in ${\rm Supp}(E)$ we deduce
$$\epsilon(D;x)=\epsilon(f^*(D);f^{-1}(x))\geq\epsilon(A;f^{-1}(x)).$$
This gives $\epsilon(D,x)\geq\epsilon(\parallel
D\parallel,x)$.  On the other hand, since $D$ is nef,
the argument in the proof of
Proposition~\ref{numerical_description} shows that we
can write
$D=A_p+\frac{1}{p}E$, with $A_p$ ample and $E$
effective, with
$x\not\in{\rm Supp}(E)$. By definition, we have
$\epsilon(\parallel D\parallel;x)
\geq\epsilon(A_p;x)$ for all $p$, and letting $p$ go to
infinity we get
$\epsilon(\parallel D\parallel;x)\geq\epsilon(D;x)$.
\end{remark}

The moving Seshadri constants measure asymptotic
separation of jets, as is the case of the usual
constants
 (see \cite{positivity}, Theorem 5.1.17). We
give now this interpretation. Recall that if $L$ is a
line bundle on a smooth variety $X$, we say that $L$
separates $s$-jets at $x\in X$ if the canonical morphism
$$H^0(X,L)\longrightarrow
H^0(X,L\otimes\OO_{X,x}/\frmm_{x}^{s+1})$$ is
surjective. Let $s(L;x)$ be the smallest $s\geq 0$ such
that
$L$ separates $s$-jets at $x$ (if there is no such
$s\geq 0$, then we put
$s(L;x)=0$).

\begin{proposition}\label{sep_jets} If $L$ is a big
line bundle on $X$, then
$$\epsilon(\parallel
L\parallel;x)=\sup_{m\to\infty}\frac{s(mL;x)}{m}
=\limsup_{m\to \infty}
\frac{s(mL;x)}{m}.$$
\end{proposition}

\begin{proof} We may clearly assume that
$x\not\in\BB(L)$, the statement being trivial otherwise.
Let $m$ be such that $x\not\in {\rm Bs}(|mL|)$, so we
have
$\pi_m : X_m\rightarrow X$ and $\pi_m^*(mL)=M_m+E_m$,
as in the definition of $\epsilon'(\parallel
L\parallel;x)$. Since $\pi_m$ is an isomorphism over a
neighborhood of $x$, it induces an isomorphism
$$mL\otimes\OO_{X,x}/\frmm_x^{s+1}\simeq
\pi_m^*(mL)\otimes
\OO_{X_m,x'}/\frmm_{x'}^{s+1},$$ where
$x'=\pi_m^{-1}(x)$. As $\pi_m$ also induces an
isomorphism
$H^0(X,mL)\simeq H^0(X_m,\pi_m^*(mL))$, we deduce
$s(mL;x)=s(\pi_m^*(mL);x')$.

On the other hand, as $x'\not\in{\rm Supp}(E_m)$,
multiplication by a local equation of $E_m$ induces an
isomorphism
$$M_m\otimes\OO_{X_m,x'}/\frmm_{x'}^{s+1}
\simeq
\pi_m^*(mL)\otimes\OO_{X_m,x'}/\frmm_{x'}^{s+1}.$$
Moreover, since $E_m$ is the fixed part of
$\pi_m^*(mL)$, we have an isomorphism
$H^0(X_m,M_m)\simeq H^0(X_m,\pi_m^*(mL))$. This gives
$s(\pi_m^*(mL);x') =s(M_m;x')$.

We show first that $\epsilon(\parallel L\parallel;x)\geq s(mL;x)/m$
for every $m$. Since $s(pmL;x)\geq p\cdot s(mL;x)$ for every $p$, we
may assume that $m$ is divisible enough, so $x \notin {\rm
Bs}(|D|)$. We take $\pi_m$ and a decomposition as above. The fact
that $\epsilon(M_m;x')\geq s(M_m;x')$ follows as in
\cite{positivity}, {\it loc. cit.}: let $C$ be an integral curve
passing through $x'$ and suppose that $M_m$ separates $s$-jets at
$x'$. We can find $F\in |M_m|$ such that ${\rm mult}_{x'}(F)\geq s$
and $C\not\subseteq F$. This gives $(F\cdot C)\geq s\cdot{\rm
mult}_xC$, hence $\epsilon(M_m;x')\geq s$. Since $\epsilon(\parallel
L\parallel;x)\geq\epsilon(M_m;x')/m$ by
Proposition~\ref{numerical_description} and since
$s(M_m;x')=s(mL;x)$, we deduce $\epsilon(\parallel L\parallel;x)
\geq s(mL;x)/m$.

In order to finish, it is enough to see that for every
$\eta>0$, we have
$s(mL;x)/m>\epsilon(\parallel L\parallel;x)-\eta$ for
some $m$. If $x \in  \B+(L)$, then the assertion
follows trivially. If $x\not\in \B+(L)$, then by
definition we can find
$f:X'\longrightarrow X$ which is an isomorphism over a
neighborhood of $x$, with
$X'$ smooth, and a decomposition $f^*(L)=A+F$, where
$A$ and $F$ are $\QQ$-divisors, with $A$ ample, $E$
effective,
$x'=f^{-1}(x)\not\in {\rm Supp}(F)$ and
$\epsilon(A;x')\geq\epsilon(\parallel
L\parallel;x)-\eta/2$.  Since $A$ is ample, it follows
from \cite{positivity}, {\it loc. cit.}, that we can
find $m$ such that $mA$ is an integral divisor and
$s(mA;x')/m\geq \epsilon(A;x')-\eta/2$. Therefore it is
enough to show that $s(mL;x)\geq s(mA;x')$. This
follows by the same arguments as before, as $\pi$ being
an isomorphism over a neighborhood of
$x$ and $x'\not\in {\rm Supp}(E)$ imply
$s(mL;x)=s(mf^*(L);x')\geq s(mA;x')$.
\end{proof}

We use Theorems~\ref{generalized_Fujita} and \ref{main}
to extend the relation between Seshadri constants and
volumes to the case of big line bundles.

\begin{proposition}\label{Seshadri_characterization} If
$D$ is a big $\QQ$-divisor on $X$ and if $x\in X$, then
$$\epsilon(\parallel D \parallel;x) = \underset{x\in
V}{\rm inf}
\frac{\vol_{X\vert V}(D)^{1/\dim(V)}}{{\rm mult}_xV},$$
where the infimum is over all positive dimensional
subvarieties $V$ containing $x$.
\end{proposition}

\begin{proof} If $x \in \B+(D)$, then $
\epsilon(\parallel D \parallel;x) = 0$ and, on the
other hand, by Theorem~\ref{main}
$\vol_{X\vert V}(D)=0$  for any irreducible component
$V$ of $\B+(D
)$ passing through $x$.

If $x \notin \B+(D)$, then any $V$ through $x$ is not
contained in $\B+(D)$, and so we can apply Theorem
\ref{generalized_Fujita}. Thus we only need to prove
that
\begin{equation}\label{eq_100}
\epsilon(\parallel D \parallel;x) = \underset{x\in
V}{\rm inf}
\frac{\parallel D^d\cdot V\parallel^{1/d}}{{\rm
mult}_xV},
\end{equation} where $d=\dim(V)$.

This however is immediate. Indeed, for each $m$
divisible enough let $\pi_m : X_m\longrightarrow X$ and
$\pi_m^*(mD)=M_m+E_m$ be as in the definition of
$\epsilon'(\parallel D\parallel;x)$, and for every $V$
denote by $\widetilde{V}_m\subseteq X_m$ the proper
transform of $V$. Since $M_m$ is nef, we have
$$\underset{x\in V}{\rm inf}
\frac{( M_m^d\cdot \widetilde{V}_m)^{1/d}}{{\rm
mult}_xV} = {\epsilon (M_m;\pi_m^{-1}(x))}.$$ It is
straightforward to deduce now equation~(\ref{eq_100})
from Proposition~\ref{numerical_description} and the
definition of
$\parallel D^d\cdot V\parallel$.
\end{proof}

The proof of the continuity of the moving Seshadri
constants is now a formal consequence of the above
results.
\begin{proof}[Proof of
Theorem~\ref{continuity_Seshadri}] Let $D$ be a
$\QQ$-divisor such that $x \notin \B+(D)$. Then the
formal concavity property of the moving Seshadri
constant  (Proposition \ref{basic1} v) gives, precisely
as in the proof of Theorem \ref{continuity} on
restricted volumes, that our function is locally
uniformly continuous around $D$. We do not repeat the
argument here.  In order to finish, it is enough to
show that if $D$ is a real class such that $x$ is in
$\B+(D)$, then
$\lim_{D'\to D}\epsilon (\parallel D'\parallel;x)=0$.
It is clear that it is enough to consider only those
$D'$ that are big. If $V$ is an irreducible component
of $\B+(D)$, then Theorem~\ref{continuity} gives
$\lim_{D'\to D}
\vol_{X\vert V}(D')=0$, so we conclude by
Proposition~\ref{Seshadri_characterization}. (This part
of the theorem is a strengthening of the main result in
\cite{nakamaye2} to the case of arbitrary real divisor
classes.)
\end{proof}

Finally, let's observe that the moving Seshadri
constant at a point controls the separation properties
of the corresponding ``adjoint'' linear series at that
point, as in the case of ample line bundles and usual
Seshadri constants (cf. \cite{demailly} or \cite{ekl}
-- the proof is essentially the same, with a slight
variation due to the initial non-positivity):
\begin{proposition}\label{effective_separation} Let $L$
be a big line bundle and assume that for some $x\in X$
$$\epsilon (\parallel L \parallel;x) > \frac{s+n}{p},$$
where $s\geq 0$ and $p > 0$ are integers, and $n= {\rm
dim}(X)$. Then the linear series $|K_X + pL|$ separates
$s$-jets at $x$.
\end{proposition}

\begin{proof} Note that the hypothesis implies that $x
\notin \B+(L)$. Let $m$ be divisible and large enough,
and let $\pi_m\colon X_m\to X$ be as in the definition
of $\epsilon'(D;x)$. If we write
$\pi_m^*(mL)=M_m+E_m$, then we have
$\epsilon(M_m;x)>\frac{m(s+n)}{p}$. For simplicity, we
identify $x$ with its inverse image in $X_m$.  We need
to prove the surjectivity of the restriction map
$$H^0 (X, \OO_X(K_X + pL))\longrightarrow H^0 (X,
\OO_X(K_X + pL)\otimes \OO_X/\frmm_x ^{s+1}).$$ Since
$K_{X_m/X}$ is supported on the exceptional locus (so,
in particular, $x$ does not lie in its support), this
is equivalent to the surjectivity of the map
$$H^0 (X_m, \OO_{X_m}(K_{X_m} +
\pi_m^*(pL)))\longrightarrow H^0 (X_m,
\OO_{X_m}(K_{X_m} + \pi_m^*(pL))\otimes
\OO_{X_m}/\frmm_x ^{s+1}).$$ This in turn is implied by
the surjectivity of the restriction map
\begin{equation}\label{surjectivity} H^0 (X_m, \OO_{X_m}
(K_{X_m} + \lceil\frac{p}{m} M_m \rceil
))\longrightarrow H^0 (X_m,
\OO_{X_m}(K_{X_m} + \lceil\frac{p}{m} M_m
\rceil)\otimes \OO_{X_m}/\frmm_x ^{s+1}),
\end{equation} since the sections on the left hand side
inject into $H^0 (X_m, \OO_{X_m}(K_{X_m} + \pi_m^*(pL)))$ by
twisting with the equation of the effective integral divisor
$\lfloor \frac{p}{m}E_m\rfloor$, and $x$ is not in the support of
this divisor. Since $M_m$ is big and nef, the  argument for
($\ref{surjectivity}$) goes as usual: consider $f: X_m' \rightarrow
X_m$ the blow-up of $X_m$ at $x$, with exceptional divisor $E$. It
is enough to prove the vanishing of $H^1 (X_m, \OO_{X_m}(K_{X_m} +
\lceil\frac{p}{m} M_m \rceil )\otimes \frmm_x^{s+1})$, which in turn
holds if
$$H^1 (X_m', f^*\OO_{X_m}(K_{X_m} + \lceil\frac{p}{m}
M_m
\rceil )\otimes \OO_{X_m'}(-(s+1)E))=0. $$ We can
rewrite this last divisor as $K_{X_m'} + f^* \lceil
\frac{p}{m} M_m\rceil - (s+n)E$, and the required
vanishing is a consequence of the Kawamata-Viehweg
Vanishing Theorem: using the lower bound on the
Seshadri constant of $M_m$ at $x$, we see that
$f^* \lceil \frac{p}{m} M_m\rceil - (s+n)E$ is the
round-up of an ample $\QQ$-divisor. (Note that we are
implicitly using here that since $M_m$ is globally
generated, the general divisor in the corresponding
linear series will avoid $x$, and so we can arrange
that $f^* \lceil \frac{p}{m} M_m\rceil =  \lceil
\frac{p}{m} f^*M_m\rceil$.)
\end{proof}

\begin{remark}
The first step in the above proof consisted in reducing to the case
when $L$ is big and nef. The rest of the proof could be
alternatively recast in the language of multiplier ideals, as
follows (see Chapter 9 in \cite{positivity} for basic facts on
multiplier ideals). Suppose that $L$ is big and nef, and fix a
rational number $t$ with $\frac{n+s}{p}<t<\epsilon(L;x)$.
Proposition~\ref{sep_jets} implies that there is $m>0$ such that
$mt$ is an integer and $mL$ separates $mt$ jets at $x$. In
particular, we can find $D\in |mL|$ whose tangent cone at $x$ is the
cone over a smooth hypersurface of degree $mt$. This implies that in
a neighborhood of $x$ we have ${\mathcal J}(qD)=\frmm_x^{s+1}$,
where $q=\frac{n+s}{mt}$. Since $(p-qm)L$ is big and nef, Nadel's
Vanishing Theorem implies that $H^1(X,\OO(K_X+pL)\otimes{\mathcal
J}(qD))=0$, hence the map
$$H^0(X,\OO(K_X+pL))\to H^0(Z,\OO(K_X+pL)\vert_Z)$$
is surjective, where $Z$ is the subscheme defined by ${\mathcal
J}(qD)$. We now deduce the assertion in the proposition using the
fact that $Z$ is defined by $\frmm_x^{s+1}$ in a neighborhood of
$x$.

\end{remark}

\smallskip

\noindent{\bf Acknowledgements}. We are grateful to Olivier Debarre
and Mihai P\u{a}un for useful discussions. We would like to thank
also the anonymous referee for several remarks and suggestions.

\end{document}